\documentclass[10pt,twoside,a4paper]{article}
\usepackage{url}
\newcommand{\liuhao}{\fontsize{8pt}{\baselineskip}\selectfont}

\usepackage{amsfonts}
\usepackage{amsmath,amssymb}
\usepackage{mathrsfs}
\usepackage{hyperref}
\usepackage{appendix}

\newtheorem{thm}{Theorem}[section]
\newtheorem{cor}[thm]{Corollary}
\newtheorem{lem}[thm]{Lemma}
\newtheorem{prop}[thm]{Proposition}

\newtheorem{rem}[thm]{Remark}

\numberwithin{equation}{section}\allowdisplaybreaks

\def\leq{\leqslant}
\def\geq{\geqslant}

\begin{document}
	\title{ {\bf \Large 
			Wellposedness and scattering for the generalized Boussinesq equation  }}
	\author
	{
		{ 
			{Jie Chen $^{1,3}$\footnote{
					Email address:
					Jiechern@163.com } }  \ \ Boling Guo$^1$ \ \	Jie Shao$^{1,2}$\footnote{Corresponding
				author;
				Email address:
				shaojiehn@foxmail.com. } } \\
		{\liuhao $^{1}$ Institute of Applied Physics and Computational Mathematics,  Beijing 100088,  P.R. China.}\\
		{\liuhao $^{2}$ Department of Mathematics, Nanjing University of Science and Technology, Nanjing 210094, P.R.China.}\\
		{\liuhao $^{3}$  School of Mathematical Sciences, Peking University,  Beijing 100871, P.R. China.}\\
		\date{}
	}
	\date{}
	\maketitle

	\begin{minipage}{13.5cm}
		\footnotesize \bf Abstract. \rm  \quad 
		In this paper, we show the local well-posedness of  the generalized Boussinesq equation(gBQ) in $L^{2}(\mathbb{R}^d), H^{1}(\mathbb{R}^d)$ and obtain the global well-posedness, finite-time blowup and small initial data scattering of gBQ in energy space $H^1(\mathbb{R}^d)$.  Moreover, we obtain the large radial initial  data scattering of defocusing case for $ d\geq 3 $ by using the method of Dodson-Murphy \cite{dodson2017new}.
		

		\vspace{10pt}
		


	\end{minipage}
	
	\section{Introduction} \label{s1.}~
	
	We study  the Cauchy problem of the generalized Boussinesq equation(gBQ)
	\begin{align} \label{ch4_1}
		\left\{\begin{array}{l}
			\partial_{t}^{2} u-\Delta u+\Delta^2 u=\beta\Delta(|u|^{\alpha-1} u), \quad (t,x)\in\mathbb{R}\times \mathbb{R}^d,\\
			u(0, x)=u_{0}(x), ~~u_{t}(0, x)=u_{1}(x),  
		\end{array}\right.
	\end{align}
	where	$ u $ is real, $ \beta=\pm1$, $\alpha>1 $. The equation has two conservation laws:
	\begin{align}
		\mathcal{E}(t): &=\frac 12\int_{\mathbb{R}^d} |(-\Delta)^{-\frac12} u_t|^2 +u^2+|\nabla u|^2+\frac{2\beta}{\alpha+1}|u|^{\alpha+1} dx=  \mathcal{E}(0), \label{ch4_est1}\\
		\mathcal{M}(t):&=\int_{\mathbb{R}^{3}}\left((-\Delta)^{-\frac{1}{2}} u_{t}\right) \nabla\left((-\Delta)^{-\frac{1}{2}} u\right) d x=\mathcal{M}(0),
	\end{align}
	and  it is
	called focusing if $ \beta=-1 $, while  it  called defocusing if $ \beta=1 $.

	The  generalized Boussinesq equation  is raised and studied in Bona et al. \cite{bona1988global} and is a generalization of the so called ``good" Boussinesq equation
	\begin{align*}
			\left\{\begin{array}{l}
				\partial_{t}^{2} u-\partial_{x}^{2} u+\partial_{x}^{4} u+\partial_{x}^{2}u^2=0, \quad (t,x)\in\mathbb{R}\times \mathbb{R}\\
			u(0, x)=u_{0}(x), ~~u_{t}(0, x)=u_{1}(x), 
		\end{array}\right.
	\end{align*}
	which models the phenomenon of nonlinear strings and is studied in Kishimoto  \cite{kishimoto2013sharp}, Mckean \cite{mckean1981boussinesq} and so on. There are several results about gBQ.  Bona et al. \cite{bona1988global} studied the local and global wellposedness and stability of solitary-wave solution in one dimension. Linares \cite{linares1993global} researched the local wellposedness  of  one dimension in $ L^2 $ and $ H^1 $ by using the $ L_p$-$L_q $ estimates. Liu \cite{liu1993instability,liu1995instability,liu1997decay} studied the local and global wellposedness, scattering, instability of solitary waves in one dimension. Cho et al. \cite{cho2007small} investigated the existence and scattering of global small amplitude
	solutions for all dimensions. Farah \cite{farah2009local} researched the local well-posedness of gBQ  for all dimensions by applying the Strichartz estimates of Schr\"odinger equation \eqref{nsl1}. Farah \cite{farah2008large,farah2012wave} studied the asymptotic behavior of solutions for  gBQ and  the spirit of the results is somewhat similar with those  in  \cite{cho2007small}.
	
	As far as we know, the results about the global wellposedness, finite-time blowup for gBQ in higher dimensions haven't been obtained yet, although  Liu \cite{liu1995instability} proved the preliminary Lemmas for all dimensions and the proof scheme of one dimension works for higher dimensions as well. What is needed is the local wellposed solution $ (u,u_t) \in H^{1} (\mathbb{R}^d )\times \dot{H}^{-1}(\mathbb{R}^d )$ that is suited for the conservation law \eqref{ch4_est1}.  Wang et al. \cite{wang2019cauchy,xmail} researched the global wellposedness, finite-time blowup, however their results depend on the damped terms $ -\alpha\Delta n_t+\gamma \Delta^2 n_{t},  \alpha\geq0, \gamma>0$  and $ -\gamma \Delta n_{t}, \gamma>0 $ respectively. The local wellposed solution 
	in Farah \cite{farah2009local}  is $   (u,(-\Delta)^{-1} u_t) \in H^{s} (\mathbb{R}^d )\times H^{s}(\mathbb{R}^d ), s\geq 0$ and their solution don't match  the results of Liu \cite{liu1995instability}. See more discussion in Appendix \ref{s5.}.

	Inspired by Gustafson, et al.\cite{gustafson2006scattering},  Kishimoto  \cite{kishimoto2013sharp}, we write \eqref{ch4_1} into
	\begin{align}
		\left(i\partial_{t}-\sqrt{(-\Delta)(1-\Delta)}\right)\left(i\partial_{t}+\sqrt{(-\Delta)(1-\Delta)}\right)u=-\beta\Delta(|u|^{\alpha-1} u).  \label{ch4_3}
	\end{align}
	and introduce
	\begin{align}
		v=u+i(-\Delta)^{-\frac12}(1-\Delta)^{-\frac12}\partial_t u=u+i\mathfrak{B}^{-1}\partial_t u, \label{trans0}
	\end{align}
	where $ \mathfrak{B}^{-1}:= (-\Delta)^{-\frac12}(1-\Delta)^{-\frac12}$.
	Then, \eqref{ch4_1} is transformed into
	\begin{align} \label{ch4_s1}
		\left\{\begin{array}{l}
			i\partial_{t} v-\mathfrak{B} v-\beta\mathfrak{M}|\textmd{Re}(v)|^{\alpha-1}\textmd{Re}(v)=0,  \\
			v(0,x)=v_0(x)=u_0+i\mathfrak{B}^{-1} u_1, 
		\end{array}\right.
	\end{align}
	where
	\begin{align}
		\mathfrak{B} &:=\sqrt{-\Delta(1-\Delta)},~
		\mathfrak{M}:=\mathscr{F}_{\xi}^{-1} \sqrt{|\xi|^{2} /\left(1+|\xi|^{2}\right)} \mathscr{F}_{x}=\sqrt{\frac{-\Delta}{1-\Delta}}, \label{ch4_trans1}\\
		u=&\frac{1}{2}\left(v+\bar{v}\right)=\textmd{Re}(v),~~
		(-\Delta)^{-\frac12}(1-\Delta)^{-\frac12}\partial_{t} u =\frac{i}{2}\left(\bar{v}-v\right)=\textmd{Im}(v). \label{ch4_trans2}
	\end{align}
	By Duhamel principle, \eqref{ch4_s1} can be written as  following integral equation,
	\begin{align}
		&v(t)= 	e^{-it \mathfrak{B}} v_{0} -i\beta\int_0^t 	e^{-i(t-\tau) \mathfrak{B}} \mathfrak{M}|\textmd{Re}(v)|^{\alpha-1}\textmd{Re}(v)d \tau, \label{ch4_inte1}
	\end{align}
	where
	\begin{align}
		e^{-it \mathfrak{B}} = \mathscr{F}_\xi^{-1} e^{-i t|\xi|\sqrt{|\xi|^2+1}} \mathscr{F}_x . \label{B}
	\end{align}

	The form  of \eqref{ch4_s1} is similar with the Cauchy problem of the nonlinear Schr\"odinger equation
	\begin{align} \label{nsl1}
		\left\{\begin{array}{l}
			i \partial_{t} z+\Delta z+\beta|z|^{\alpha-1} z=0, \quad t \geq 0, x \in \mathbb{R}^{d} \\
			z(0,x)=z_0(x), 
		\end{array}\right.
	\end{align}
	where $\alpha>1, \beta=\pm1$.  The nonlinear Schr\"odinger equation \eqref{nsl1} is one of the most famous dispersive equations and there are numerous results and many mature theories about it.   We refer to Cazenave \cite{cazenave2020overview} for a brief introduction of \eqref{nsl1}.	Naturally, we shall try to derive results that are analogous to those of \eqref{nsl1}.  The result in Gustafson et al. \cite{gustafson2006scattering} has given  the Strichartz type estimates for operator \eqref{B}, which is presented in Section \ref{s2.}, and this give us the start point and an effective tool  to research the results in this paper. The local wellposedness is  directly obtained by following the classic results studying the Schr\"odinger equation and the global wellposednes, finite-time blowup thus follows by basing on the analysis of Liu \cite{liu1995instability} and Wang et al. \cite{xmail}.
	
	There are rare results about the scattering of gBQ. Liu \cite{liu1997decay} researched the small initial data scattering in one dimension. Cho et al. \cite{cho2007small} researched the small initial data scattering in higher dimensions. Mu\~noz et al. \cite{munoz2018scattering} studied the small initial data scattering in one dimension  as well, and their method is much different from the classic ones studying the Schr\"odinger equation, wave equation etc.  Actually, their way is more likely as those investigating stability of soliton solutions. It is not difficult to study small initial data scattering by using the transformation \eqref{trans0} and  the Strichartz type estimates in Lemma \ref{le2.1}. However, for the large initial data scattering, it is quite nontrivial. The $ \Delta $ in the nonlinear term $ \Delta(|u|^{\alpha-1} u) $ makes it difficult to calculate the classic Morawetz type estimate  studying the Schr\"odinger equation. Specifically, if we research on \eqref{ch4_s1}, the operator  $ 	\mathfrak{M} $  before the nonlinear
	term $|\textmd{Re}(v)|^{\alpha-1}\textmd{Re}(v) $ makes it hard to construct   Morawetz type estimate  perfect as those for the Schr\"odinger equation. Inspired by Dodson et al. \cite{dodson2017new}, we calculate the Morawetz-virial type estimate related to \eqref{ch4_1}.
	We remark that the calculation relies heavily  on the radial Sobolev embedding inequality and the non-radial large data scattering is still open.

	The main results in this paper are as follows.
	\begin{thm} \label{th1}  Local wellposedness
		
		The initial value problem \eqref{ch4_1} with $(u, (-\Delta)^{-\frac12}u_t)$ in $H^{s} (\mathbb{R}^{d}) \times H^{s-1} (\mathbb{R}^{d})$ is locally well-posed  for
		\begin{align}  
			\begin{array}{ll}
				s=0, &  \text { if } 1<\alpha\leq1+4 /d, \\
				s = 1 , & \text { if } 
				1<\alpha<\infty,  d=1,2
				\text{  or  } 1<\alpha\leq\frac{d+2}{d-2},   d\geq 3.
			\end{array}
		\end{align}	
		
	\end{thm}
	
	Before introducing Theorem \ref{th2}, we bring in some definitions used in Liu \cite{liu1995instability} and the notations may be slightly different from those in \cite{liu1995instability}.
	\begin{equation} \label{nota}
		\begin{aligned}
			E(u):=\frac12\left\|u\right\|_{H^1(\mathbb{R}^d)}^2-\frac{1}{\alpha+1}\left\|u\right\|_{L^{\alpha+1}(\mathbb{R}^d)}^{\alpha+1}, \\
			R(u):=\left\|u\right\|_{H^1(\mathbb{R}^d)}^2-\left\|u\right\|_{L^{\alpha+1}(\mathbb{R}^d)}^{\alpha+1},\\
			\eta_d:=\min\left\{	E(u)\left.\right| 0\neq u\in H^1(\mathbb{R}^d),R(u)=0\right\},
		\end{aligned}
	\end{equation}
	$ \varphi_d $ is the ground state solution of
	\begin{align}
		-\Delta \varphi_d+\varphi_d-|\varphi_d|^{\alpha-1}\varphi_d=0,\quad x\in \mathbb{R}^d. \label{ground}
	\end{align}
	Define
	\begin{equation} \label{nota1}
			C_\ast=C_\ast(\alpha,d):=	C_{\alpha,d}
	\end{equation}
	as the best constant for Sobolev inequality $ \| u\|_{L^{\alpha+1}(\mathbb{R}^d)}  \leq  	C_{\alpha,d}	\| u\|_{H^1(\mathbb{R}^d)} $ and it is determined in Proposition \ref{propg}.
	If $ \beta=-1 $, then $ 	\mathcal{E}(t)=E(u)+\frac12\|u_t\|_{\dot{H}^{-1}(\mathbb{R}^d)}^2. $ 
	\begin{thm}  \label{th2} Global wellposedness and finite-time blowup\\	
		$ (i) $  	Assume $ 	\mathcal{E}(0) $ is a finite positive constant and $ \alpha $ satisfies
		\begin{align} \label{ch4_alpha0}
			\left\{\begin{array}{ll}
				1<\alpha<\infty, & \text { if } d=1,2,\\
				1<\alpha<\frac{d+2}{d-2}, & \text {if } d\geq 3.			
			\end{array}\right.
		\end{align} 			
		Then, the initial value problem \eqref{ch4_1} with $(u, (-\Delta)^{-\frac12}u_t)$ in $H^{s} (\mathbb{R}^{d}) \times H^{s-1} (\mathbb{R}^{d})$ is globally wellposed for $\beta = 1$ or 
		\begin{align}  \label{assum1}
			~~\beta=-1,~\mathcal{E}(0)<\frac{\alpha-1}{2(\alpha+1)}C_\ast^{-\frac{2(\alpha+1)}{\alpha-1}},\quad  \|u_0\|_{H^1}<C_\ast^{-\frac{\alpha+1}{\alpha-1}}. 		
		\end{align} 	
		\noindent{$(ii)$} 	Assume	$ \alpha $ satisfies \eqref{ch4_alpha0}, $ \beta=-1 $,   $ (u_0,u_1)\in H^{1}(\mathbb{R}^{d})\times \dot{H}^{-1}(\mathbb{R}^{d}),  u_0\in \dot{H}^{-1}(\mathbb{R}^{d})$ and 
		\begin{align} \label{assum2}
			\mathcal{E}(0)<\frac{\alpha-1}{2(\alpha+1)}C_\ast^{-\frac{2(\alpha+1)}{\alpha-1}},~~ \|u_0\|_{H^1}>C_\ast^{-\frac{\alpha+1}{\alpha-1}},
		\end{align}  	
		then  the  local wellposed solution of \eqref{ch4_1}  $ u\in H^{1}(\mathbb{R}^{d})$ blows up in some finite time, i.e. the maximum lifespan $ T_{max}<\infty  $,  and $ u  $ satisfies
		\begin{align*}
			\lim\limits_{t\rightarrow T_{max}^-} \left\|u(t)\right\|_{H^1(\mathbb{R}^d)}=	\lim\limits_{t\rightarrow T_{max}^-} \left\|u(t)\right\|_{L^{\alpha+1}(\mathbb{R}^d)}=\infty.
		\end{align*} 	
		\noindent{$(iii)$} 	Let $ \varphi_d \in H^1 (\mathbb{R}^d) $ be the ground state of \eqref{ground} and assume  $ u_1=0, \beta =-1 $. For any $ \delta>0 $ , there exists initial data $ u_0 \in  H^1 (\mathbb{R}^d) $ with $ \|u_0-\varphi_d\|_{H^1}<\delta $, such that  the local 
		solution of \eqref{ch4_1}  
		in $H^{1}(\mathbb{R}^{d})\times \dot{H}^{-1}(\mathbb{R}^{d}) $  satisfies  	
		\begin{align*}
			\lim\limits_{t\rightarrow T^-} \left\|u(t)\right\|_{H^1}=\infty
		\end{align*} 
		for some $ 0<T<\infty $.
		
	\end{thm}

	\begin{thm} \label{th3} Small initial data scattering~
		
		Suppose $ (u_0,u_1)\in H^{1}(\mathbb{R}^{d})\times \dot{H}^{-1}(\mathbb{R}^{d}) $ small enough and $ \alpha $ satisfying
		\begin{align} \label{ch4_alpha1}
			\left\{\begin{array}{ll}
				4/d+1\leq\alpha<\infty, & \text { if } d=1,2,\\
				4/d+1\leq\alpha\leq\frac{d+2}{d-2}, & \text {if } d\geq 3,			
			\end{array}\right.
		\end{align} 
		then the   $ H^1 $ solution of \eqref{ch4_1} $ (u,u_t) \in H^{1}(\mathbb{R}^{d})\times \dot{H}^{-1}(\mathbb{R}^{d}) $ scatters. More precisely, there exists unique linear solution of \eqref{ch4_1} $ (u_0^\pm,u_1^\pm)\in H^{1}(\mathbb{R}^{d})\times \dot{H}^{-1}(\mathbb{R}^{d}) $  such that
		\begin{align}
			\lim _{t \rightarrow \pm \infty}\left\|\left(u(t), u_{t}(t)\right)-\left(u_0^{\pm}, u_1^{\pm}\right)\right\|_{\mathcal{H}^{1}}=0, \label{sc1}
		\end{align}
		where $\mathcal{H}^{1}:=H^1\times\dot{H}^{-1}$,
		$\left\|\left(w_1,w_2\right)\right\|_{\mathcal{H}^{1}}:=(\|w_1\|_{H^{1}}^{2}+\|w_2\|_{\dot{H}^{-1}}^{2})^{\frac{1}{2}}=\left\|w_1+i \mathfrak{B}^{-1} w_2\right\|_{H^{1}}$.
	\end{thm}
	
	\begin{thm} \label{th4} Large radial initial data scattering  of defocusing case ~
		
		Suppose $ (u_0, u_1) $ is radial symmetric, $ \beta=1 $,  $ 4/d+1\leq\alpha<(d+2)/(d-2) $ and $ d\geq 3 $,   then the $ H^1 $ solution of \eqref{ch4_1} scatters.
	\end{thm}
	
	The rest of the paper is organized as follows. We research the local and global wellposedness, finite-time blowup in Section \ref{s2.} and derive the small initial data scattering  in Section \ref{s3.}. The large initial data scattering for  radial defocusing case is obtained in Section \ref{s4.}. Some more discussion about the properties of the Boussinesq operator is given in Section \ref{s5.}.
	
	We end up with notations used in this paper. $\mathscr{F}_{x}u$ and $\widehat{u}$ denote the Fourier transform
	of the function $u(t,x)$ with respect to $x$, and $\mathscr{F}_{\xi}^{-1}u, u^{\vee}(x)$ the inverse Fourier transformation to $\xi$ :
	$$\mathscr{F}_{x}u=\widehat{u}(\xi):=\int_{\mathbb{R}^{d}} u(x) e^{-2 \pi \mathrm{i} x\cdot \xi} dx, ~~\mathscr{F}_{\xi}^{-1}u=u^{\vee}(x):=\widehat{u}(-x).$$
	Define  $ |x|:=\sqrt{\Sigma_i x_i^2} $ and
	$ D^\eta n:= \mathscr{F}^{-1}|\xi|^\eta (\mathscr{F} n )(\xi),
	\eta \in \mathbb{R},  \xi \in \mathbb{R}^d, $ then $ D=(-\Delta)^{\frac 12}$.
	Denote  $ \left\langle\cdot\right\rangle:=\sqrt{1+|\cdot|^2},~\left\langle D\right\rangle:=\sqrt{1+D^2},~\left\langle\nabla\right\rangle$  $:=1+\nabla $, then we have $  \left\langle D \right\rangle \sim\left\langle\nabla\right\rangle$. Write $A \leq C B, A \leq C(a) B$ as $A \lesssim B , A \lesssim_a B  $ respectively, where  $  C , C(a) $  are positive constants and $ C(a) $ depends on $ a $. The Rezis operator  is denoted as $  R_{j}:=-\mathscr{F}_\xi^{-1} \frac{i\xi_{j}}{|\xi|} \mathscr{F}_x $. The Lebesgue, Sobolev and Besov norms are standard
	\begin{align*}
		\|u\|_{ L^{p}(\mathbb{R}^d)}=\left(\int_{\mathbb{R}^d}|u(x, t)|^{p} d x\right)^{\frac{1}{p}}, \|u\|_{L^q([0,t]; L^{p}(\mathbb{R}^d))}=\left(\int_{0}^{t}\|u\|_{ L^{p}(\mathbb{R}^d)} ^q d t\right)^{\frac{1}{q}}&,\\
		\|u\|_{W^{\gamma,q}(\mathbb{R}^d)}=\left\|\langle D\rangle^{\gamma} u\right\|_{L^{q}(\mathbb{R}^d)},~ \|u\|_{H^{\gamma}(\mathbb{R}^d)}=\left\|\langle D\rangle^{\gamma} u\right\|_{L^{2}(\mathbb{R}^d)},~ \|u\|_{\dot{H}^{\gamma}(\mathbb{R}^d)}=\left\| D^{\gamma} u\right\|_{L^{2}(\mathbb{R}^d)}&,\\
		\|u\|_{\dot{B}_{p, r}^{s}}=\left(\sum_{j \in \mathbb{Z}} 2^{r j s}\left\|\dot{\Delta}_{j} u\right\|_{L^{p}}^{r}\right)^{\frac{1}{r}},~\dot{\Delta}_{j} u=\varphi\left(2^{-j} D\right) u, ~ \sum_{j \in \mathbb{Z}} \varphi\left(2^{-j} \xi\right)=1, ~\forall \xi \in \mathbb{R}^{d} \backslash\{0\},&
	\end{align*}
	and $ \varphi $ is the Littlewood-Paley function, which is radial symmetry and supported in the  annulus $ \left\{\xi \in \mathbb{R}^{d} ~|~ 3 / 4 \leq|\xi| \leq 8 / 3\right\} $. If there is no special explanation, the norm without domain subscript is always acting on the whole space $ \mathbb{R}^d $. 
	
	\section{Local and global wellposedness, finite-time blowup} \label{s2.}
	\subsection{Local  wellposedness} 
	\quad
	
	The Strichartz type estimates for operator \eqref{B} in Gustafson et al. \cite{gustafson2006scattering} are as follows.
	\begin{lem}  \label{le2.1}
		Let $d \geq 1 $ be integer. \\
		$ (i) $ For $2 \leq q \leq \infty$, we have
		
		\begin{align}
			\left\|e^{- it \mathfrak{B}} \varphi\right\|_{\dot{B}_{q,2}^{0}} \lesssim |t|^{-d \sigma}\left\|\mathfrak{M}^{(d-2) \sigma} \varphi\right\|_{\dot{B}_{q^{\prime},2}^{0}} \label{stype1}
		\end{align}
		where $\frac1q+\frac{1}{q^{\prime}}=1$, $\sigma=1 / 2-1 / q$.
		\\
		$ 	(ii) $ Suppose $j=1,2,$  $2 \leq p_{j}, q_{j} \leq \infty, 2 / p_{j}+d / q_{j}=d / 2$, $s_{j}=\frac{d-2}{2}(1 / 2-1 / q_{j}),$ but $\left(p_{j}, q_{j}\right) \neq(2, \infty) $, then
		\begin{align}
			\left\|e^{-i t \mathfrak{B}} \varphi\right\|_{L^{p_{1}} \dot{B}_{q_{1},2}^{0}} &\lesssim\left\|\mathfrak{M}^{s_{1}} \varphi\right\|_{L^{2}}, \label{stype2}\\
			\left\|\int_{-\infty}^{t} e^{-i(t-s)  \mathfrak{B}} f(s) d s\right\|_{L^{p_{1}} \dot{B}_{q_{1},2}^{0}}& \lesssim\left\| \mathfrak{M}^{s_{1}+s_{2}} f\right\|_{L^{p_{2}^{\prime}} \dot{B}_{q_{2}^{\prime},2}^{0}},\label{stype3}
		\end{align}
		where $ \mathfrak{M} $ is defined as \eqref{ch4_trans1}.
	\end{lem}
	
	Analogous to the Strichartz estimates for \eqref{nsl1}, we call  $ (p,q) $ is  admissible, if it satisfies 
	\begin{equation} \label{sadmiss}	
		\begin{aligned} 
			&\qquad	\frac{2}{p}=d\left(\frac{1}{2}-\frac{1}{q}\right), \\	
			&	\left\{\begin{array}{ll}
				2 \leq q \leq \infty, &  \text { if } d=1, \\
				2 \leq q<\infty, & \text { if } d=2, \\
				2 \leq q \leq2 d/{(d-2)}, & \text { if } d \geq 3.
			\end{array}\right.
		\end{aligned}
	\end{equation}	
	By the definition of $ s_i $ and the relationship between Sobolev space and Besov space, one can follow \eqref{stype2}, \eqref{stype3} to get
	\begin{align}
		\left\|e^{-i t \mathfrak{B}} \varphi\right\|_{L^{p_{1}} L^{q_1}} &\lesssim\left\|  \mathfrak{M}^{s_{1}}\varphi\right\|_{L^{2}}, \label{stype4}\\
		\left\|\int_{-\infty}^{t} e^{-i(t-s)  \mathfrak{B}} f(s) d s\right\|_{L^{p_{1}} L^{q_1}}&\lesssim\left\|\int_{-\infty}^{t} e^{-i(t-s)  \mathfrak{B}} \mathfrak{M}^{s_{1}+s_{2}}f(s) d s\right\|_{L^{p_{1}} \dot{B}_{q_{1},2}^{0}}\nonumber\\
		& \lesssim\left\|  \mathfrak{M}^{s_{1}+s_{2}} f\right\|_{L^{p_{2}^{\prime}} \dot{B}_{q_{2}^{\prime},2}^{0}}\lesssim\left\| \mathfrak{M}^{s_{1}+s_{2}} f\right\|_{L^{p_{2}^{\prime}}L^{q_{2}^{\prime}}  }, \label{stype5}
	\end{align}
	$ 2\leq q_i <\infty, 1<q_i^\prime\leq 2, q_i=1,2$. One can see \eqref{stype4}, \eqref{stype5} and  the Strichartz estimates  are very similar in form and 	\eqref{stype4}, \eqref{stype5}  are essentially better. Indeed,  formally, there is 
	\begin{equation*}
		\left|	\mathscr{F}\left(\mathfrak{M}^{s} g(x) \right) \right| \leq 	\left|	\widehat{g} \right|, ~s>0
	\end{equation*}
and for the multiplier of operator $  \mathfrak{M}^{s} , s>0 $, namely $(\frac{|\xi|^2}{1+|\xi|^2})^{\frac s2} $, it is easy to calculate that 
\begin{align*}
	\left|\partial_{\xi}^{\alpha} (\dfrac{|\xi|^2}{1+|\xi|^2})^{\frac s2}\right| \leq A|\xi|^{-|\alpha|}
\end{align*}
for $ |\alpha| \leq\left[\frac{d}{2}\right]+1 $ and some constant $ A $, thus  the H\"ormander–Mihlin Multiplier Theorem shows $(\frac{|\xi|^2}{1+|\xi|^2})^{\frac s2} $ is a $ L^p,  1<p<\infty $ multiplier and $ \|\mathfrak{M}^{s}f\|_{L^p}\lesssim \|f\|_{L^p} $.  
	On the other hand, for the high frequency, there is
	\begin{equation*}
		\left|	\mathscr{F}\left(\mathfrak{M}^{s}g(x) \right) \right|  \sim	\left|	\widehat{g} \right|, 
	\end{equation*} 
and it seems that we can't expect  results for \eqref{ch4_s1} that are better than those of \eqref{nsl1} within the standard method for discussing local wellposdness for \eqref{nsl1}. 
	Thus, the classification for  \eqref{nsl1} according to exponent $\alpha$ of nonlinear term, such as  mass critical case $\alpha=1+4/d$ and  energy critical case $\alpha=1+4/(d-2), d \geq 3$,  may be still used to
	sort through   \eqref{ch4_s1} or \eqref{ch4_1}. 
	
	The local wellposedness results for \eqref{ch4_s1} is presented as below.	
	\begin{lem} \label{le2}
		$1<\alpha<1+4 /d$ case in $L^{2}\left(\mathbb{R}^{d}\right)$ 
		
		If $1<\alpha<1+4 /d,$ then for every $v_{0} \in L^{2}\left(\mathbb{R}^{d}\right)$,
		there exist $T=T\left(\left\|v_{0}\right\|_{L^2}, d,  \alpha\right)>0$ and a unique solution $v$ of the integral equation \eqref{ch4_inte1} in the time interval $[-T, T]$ with
		$$
		v \in C\left([-T, T]; L^{2}\left(\mathbb{R}^{d}\right)\right) \cap L^{r}\left([-T, T]; L^{\alpha+1}\left(\mathbb{R}^{d}\right)\right),
		$$
		where $r=\frac{4(\alpha+1)}{ d(\alpha-1)}$.
		Moreover, for every $T^{\prime}<T$ there exists a neighborhood $V$ of $v_{0}$ in $L^{2}\left(\mathbb{R}^{d}\right)$ such that
		$$
		\mathbb{F}: V \mapsto C\left(\left[-T^{\prime}, T^{\prime}\right]; L^{2}\left(\mathbb{R}^{d}\right)\right) \cap L^{r}\left(\left[-T^{\prime}, T^{\prime}\right]; L^{\alpha+1}\left(\mathbb{R}^{d}\right)\right), \quad \tilde{v}_{0} \mapsto \tilde{v}(t)
		$$
		is Lipschitz.
	\end{lem}
	
	\begin{lem} \label{le3}
		Critical case, $\alpha=1+4 / d$ in $L^{2}\left(\mathbb{R}^{d}\right)$
		
		If $\alpha=1+4 /d,$ then for every  $v_{0} \in L^{2}\left(\mathbb{R}^{d}\right)$ there exist $T=T\left(v_{0},  \alpha\right)>0$ and a unique solution $v$ of the integral equation \eqref{ch4_inte1} in the time interval $[-T, T]$ with
		$$
		v \in C\left([-T, T]; L^{2}\left(\mathbb{R}^{d}\right)\right) \cap L^{\sigma}\left([-T, T ]; L^{\sigma}\left(\mathbb{R}^{d}\right)\right),
		$$
		where $\sigma=2+4 / d$. Moreover, for every $T^{\prime}<T$ there exists a neighborhood $V$ of $v_{0}$ in $L^{2}\left(\mathbb{R}^{d}\right)$ such that
		$$
		\mathbb{F}: V \mapsto C\left([-T^{\prime}, T^{\prime}]; L^{2}\left(\mathbb{R}^{d}\right)\right) \cap L^{\sigma}\left(\left[-T^{\prime}, T^{\prime}\right]; L^{\sigma}\left(\mathbb{R}^{d}\right)\right), \quad \tilde{v}_{0} \mapsto \tilde{v}(t)
		$$
		is Lipschitz.
	\end{lem}	
	\begin{lem} \label{le4} Local theory in $H^{1}\left(\mathbb{R}^{d}\right)$
		
		If $\alpha$ 	satisfying
		\begin{align} \label{alpha1}
			\left\{\begin{array}{ll}
				1<\alpha<\frac{d+2}{d-2}, & \text { if } d\geq 3 \\
				1<\alpha<\infty, & \text { if } d=1,2,
			\end{array}\right.
		\end{align} 
		then for every $v_{0} \in H^{1}\left(\mathbb{R}^{d}\right)$, there exist $T=T\left(\left\|v_{0}\right\|_{H^1}, d,  \alpha\right)>0$ and a unique solution v of the integral equation \eqref{ch4_inte1} in the time interval $[-T, T]$ with
		$$
		v \in C\left([-T, T]; H^{1}\left(\mathbb{R}^{d}\right)\right) \cap L^{p}\left([-T, T]; W^{1,q}\left(\mathbb{R}^{d}\right)\right),
		$$
		where $(p, q)=\left(\frac{4(\alpha+1)}{d(\alpha-1)},\alpha+1\right)$ for $d \geq 3,$ and $(p, q)$ satisfies \eqref{sadmiss}
		for $d=1,2$.
				
		Moreover, v depends continuously on $ v_0(x) $ as follows. Let 
		$\left(v_{0,n}\right)_{n \geq 1} \subset$ $H^{1}\left(\mathbb{R}^{d}\right)$ such that $ v_{0,n} \rightarrow  v_0 $ in $H^{1}\left(\mathbb{R}^{d}\right)$ as $n \rightarrow \infty$, and let $v_{n}$ be the maximal solution of  \eqref{ch4_inte1} corresponding to the initial value $ v_{0,n}$, then there exists $0<T^{\prime}<T$ depending on $\|v_0\|_{H^{1}}$ such that $v_{n}$ is defined on $[-T^{\prime}, T^{\prime}]$ for $n$ large enough and $v_{n} \rightarrow v$ in $C\left([-T^{\prime}, T^{\prime}], H^{1}\left(\mathbb{R}^{d}\right)\right)$ as $n \rightarrow \infty$.
	\end{lem}
	
		

	\begin{lem} \label{le6}  
		Critical case, 
				$\alpha=(d+2) /(d-2), d\geq 3,$ in $H^{1}\left(\mathbb{R}^{d}\right)$
		
				Let $d\geq 3$,
			$\alpha=(d+2) /(d-2) $, $p=\alpha+1=2 d /(d-2)$ and $ q=\frac{2d(\alpha+1)}{d(\alpha+1)-4}=2 d^{2} /\left(d^{2}-2 d+4\right)$, then for every $v_{0} \in H^{1}\left(\mathbb{R}^{d}\right)$, there exist $T=T\left(v_{0}, d,  \alpha\right)>0$ and a unique solution $ v $ to \eqref{ch4_inte1} satisfying
			\begin{align*}
				v \in C\left(\left[-T, T\right]; H^{1}(\mathbb{R}^{d})\right) \cap L_{\mathrm{loc}}^{p}\left(\left[-T, T\right]; W^{1,q}(\mathbb{R}^{d})\right).
			\end{align*}
	
		Moreover, v depends continuously on $ v_0(x) $  as follows. Let $\left(v_{0,n}\right)_{n \geq 1} \subset$ $H^{1}\left(\mathbb{R}^{d}\right)$ such that $ v_{0,n} \rightarrow  v_0 $ in $H^{1}\left(\mathbb{R}^{d}\right)$ as $n \rightarrow \infty$, and let $v_{n}$ be the maximal solution of  \eqref{ch4_inte1} corresponding to the initial value $ v_{0,n}$, then exists $0<T^{\prime}<T$ depending on $v_0$ such that $v_{n}$ is defined on $[-T^{\prime}, T^{\prime}]$ for $n$ large enough and $v_{n} \rightarrow  v$ in $L^\infty\left([-T^{\prime}, T^{\prime}], H^{1}\left(\mathbb{R}^{d}\right)\right)$, as  $n \rightarrow \infty$.
		
	\end{lem}
	\begin{flushright}
		$\lozenge$
	\end{flushright}
	
	The results list can go on. Due to the similarity between \eqref{ch4_s1} and \eqref{nsl1}, we will only prove  Lemma \ref{le6} and refer the proof of  Lemma \ref{le2}-\ref{le3} to  Theorem 5.2, 5.3 in Linares et al. \cite{linares2014introduction} respectively and Lemma \ref{le4} to Theorem 4.4.1 in Cazenave \cite{2003semilinear}.  The exponent $ \alpha=(d+2) /(d-2) $  is  critical  for Schr\"odinger equation \eqref{nsl1} in energy space $  \dot{H}^{1}$, as then
	\eqref{nsl1} is invariant in $ \dot{H}^{1} $ under the scaling  $
	z_{\lambda}(t, x)=\lambda^{-\frac{2}{\alpha-1}} z\left(\lambda^{-2} t, \lambda^{-1} x\right)
	$, i.e.
	\begin{align*}
		\left\|z_{\lambda}(0,x)\right\|_{\dot{H}^{1}}=\lambda^{\frac{d}{2}-\frac{2}{\alpha-1}-1}\|z(0,x)\|_{\dot{H}^{1}}=\|z(0,x)\|_{\dot{H}^{1}}.
	\end{align*}
Both of  \cite{2003semilinear} and \cite{linares2014introduction}  have given the  local wellposedness for \eqref{nsl1} of energy critical case in $ H^1(\mathbb{R}^d) $. The  method of Theorem 4.5.1 in \cite{2003semilinear} obtaining the results  is  an iteration in Strichartz spaces and those of Theorem 5.5 in \cite{linares2014introduction} is  Kato's method. However, as clarified in Tao et al. \cite{2005Stability}, it seems if the dimension  $ d $ is large enough, we can't get the Lipschitz continuity of solution $ z(t,x) $ for \eqref{nsl1} of energy critical case in $ H^1(\mathbb{R}^d)$ that are stated in Theorem 5.5 of \cite{linares2014introduction}. Here,  the proof method for Lemma \ref{le6} is a combination of those in \cite{2003semilinear} and \cite{linares2014introduction}. We will use Kato's method and  prove the continuity for the solution $ v $ in the sense  as those in \cite{2003semilinear}.

	We  recall the Kato-Ponce type inequalities to deal with the nonlinear terms. 
	\begin{lem} \cite{grafakos2014kato}
		Let $\frac{1}{2}<r<\infty, \quad 1<p_{1}, p_{2}, q_{1}, q_{2} \leq \infty$ satisfy $\frac{1}{r}=\frac{1}{p_{1}}+\frac{1}{q_{1}}=\frac{1}{p_{2}}+\frac{1}{q_{2}} .$
		Given $s>\max \left(0, \frac{d}{r}-d\right)$ or $s \in 2 \mathbb{N}$, there exists $C=C\left(d, s, r, p_{1}, q_{1}, p_{2}, q_{2}\right)<\infty$
		such that for all $f, g \in \mathcal{S}\left(\mathbb{R}^{d}\right)$, we have
		\begin{align}
			\left\|D^{s}(f g)\right\|_{L^{r}\left(\mathbb{R}^{d}\right)} \lesssim &\left\|D^{s} f\right\|_{L^{p_{1}}\left(\mathbb{R}^{d}\right)}\|g\|_{L^{q_{1}}\left(\mathbb{R}^{d}\right)}+\|f\|_{L^{p_{2}\left(\mathbb{R}^{d}\right)}}\left\|D^{s} g\right\|_{L^{q_{2}}\left(\mathbb{R}^{d}\right)},  \label{nonlinear1}\\
			\left\|\left\langle D \right\rangle^s(f g)\right\|_{L^{r}\left(\mathbb{R}^{d}\right)} & \lesssim\|f\|_{L^{p_{1}}\left(\mathbb{R}^{d}\right)}\left\|\left\langle D \right\rangle^s g\right\|_{L^{q_{1}}\left(\mathbb{R}^{d}\right)}+\left\|\left\langle D \right\rangle^s f\right\|_{L^{p_{2}}\left(\mathbb{R}^{d}\right)}\|g\|_{L^{q_{2}}\left(\mathbb{R}^{d}\right)}. \label{nonlinear2}
		\end{align}
	\end{lem}
We also need following Lemma to deal with the initial data term. 
\begin{prop} \label{propini}
 Let $(p, q)$ satisfy \eqref{sadmiss}. Given $u_{0} \in L^{2}\left(\mathbb{R}^{d}\right)$ and $\varepsilon>0$, there exist $\delta>0$ and $T>0$ such that if $\left\|v_{0}-u_{0}\right\|_{2}<\delta$, then
\begin{align} \label{propini1}
		\left(\int_{0}^{T}\left\|e^{-it \mathfrak{B}}v_{0}\right\|_{q}^{p} d t\right)^{1 / p}<\varepsilon.
\end{align}
\end{prop}
Noticing that $ \text { that }\left\{e^{-it \mathfrak{B}}\right\} \text { defines a unitary group in } H^{s}\left(\mathbb{R}^{d}\right) $, Proposition \ref{propini}  is essentially the same with Proposition 5.1 of \cite{linares2014introduction}.
	
	\noindent{\bf{Proof of Lemma \ref{le6}:}}
	
Without loss of generality, we set $  t\geq 0 $. Define space as
	\begin{align*}
		X_T:=&\left.\left\{v \in L^\infty\left([0, T]; H^{1}\left(\mathbb{R}^{d}\right)\right)\cap L^{p}\left([0,T]; W^{1,q}\right)\right.\right|\\
		&~\sup _{\left[0, T\right]}\left\|v(t)-e^{-it \mathfrak{B}} v_{0}\right\|_{H^1}+\|v\|_{L^{p}\left([0,T]; W^{1,q}\right)} \leq M\},
	\end{align*}
	and distance norm of $ X_T $ as	$$
	\mathfrak{d}(v, w)=\sup _{\left[0, T\right]}\left\|v-w\right\|_{L^2}+\|v-w\|_{L^{p}\left([0,T]; L^{q}\right)}.
	$$
	By  Theorem 1.2.5 of Cazenave \cite{2003semilinear} and following  corresponding discussion in  Theorem 4.4.1 of \cite{2003semilinear}, we know $ (X_T,\mathfrak{d}) $ is a complete metric space. Thus,
	one only needs to prove
	\begin{align}
		\mathscr{T}: v \mapsto e^{-it \mathfrak{B}} v_{0} -i\beta\int_0^t 	e^{-i(t-\tau) \mathfrak{B}} \mathfrak{M}|\textmd{Re}(v)|^{\alpha-1}\textmd{Re}(v)d \tau, \label{ch4_inte2}
	\end{align}
	is a contraction map on $ (X_T,\mathfrak{d}) $.
	
	As $ v_0\in H^{1}(\mathbb{R}^{d}) $,   we can let the $ u_0,v_0 $ in Proposition \ref{propini} be $ v_0=u_0 $ and find small enough $ T>0 $ such that
	\begin{align*}
			\left(\int_{0}^{T}\left\|e^{-it \mathfrak{B}} v_{0}\right\|_{W^{1,q}}^p d t\right)^{1 / p}<\varepsilon,
	\end{align*}
	for any given small $ \varepsilon>0 $, where
$p=\alpha+1=2 d /(d-2)$ and  $ q=\frac{2d(\alpha+1)}{d(\alpha+1)-4}=2 d^{2} /\left(d^{2}-2 d+4\right)$ is admissible. Let $ r_1=d(\alpha+1)(\alpha-1)/4= 2 d^{2} /\left(d-2\right)^2$, then the Sobolev embedding inequality
	\begin{align*}
		\frac{1}{q^\prime}=\frac1q+\frac{\alpha-1}{r_1},~~W^{{1},q}(\mathbb{R}^{d})\hookrightarrow L^{r_1}(\mathbb{R}^{d}),
	\end{align*} 
 by nonlinear estimate \eqref{nonlinear2}, Strichartz  estimates \eqref{stype4},  \eqref{stype5} and H\"older inequality show
 \begin{align} \label{loc}
 		\left\|\mathscr{T}	(v)\right\|_{L^{p}([0,T];W^{1,q})}&\leq C
 		\left\|e^{-it \mathfrak{B}}v_0\right\|_{L^{p}([0,T];W^{1,q})}+C	\left\|\mathfrak{M}^{1+2s_1}|\textmd{Re}(v)|^{\alpha-1}\textmd{Re}(v)\right\|_{L^{p^\prime}([0,T];W^{{1},{q^\prime}})}\nonumber\\
 		&	\leq C
 	\varepsilon+ C\left\|\left\||v|^{\alpha-1}\right\|_{L^{\frac{d(\alpha+1)}{4}}}	\left\|v\right\|_{W^{{1},{q}}}\right\|_{L^{p^\prime}[0,T]}\nonumber\\
 		&\leq C
 	\varepsilon+ C	\left\|v\right\|_{L^{p}([0,T];W^{{1},{q}})}^\alpha,
 \end{align}
where $ s_{1}=\frac{d-2}{2}(1 / 2-1 / q) $ is determined as those in  Lemma \ref{le2.1}. 
 By Bootstrap argument and the smallness of $ \varepsilon $, we know $ \left\|v\right\|_{L^{p}([0,T];W^{{1},{q}})} \leq C(\varepsilon) $ is small and can set $ M $  small enough such that 
 \begin{align} \label{loc0}
 	C 
 	M^{\alpha-1} \leq \frac12, C \varepsilon+CM^{\alpha} \leq M.
 \end{align}
 Therefore, 
 \begin{align*}
	\left\|\mathscr{T}	(v)\right\|_{L^{p}([0,T];W^{1,q})} \leq C
		\varepsilon+ C M^\alpha \leq M
 \end{align*}
and similar discussion with \eqref{loc} gives
	\begin{align}
\begin{aligned} \label{loc1}
	\sup _{\left[0, T\right]}\left\|\mathscr{T}	(v)-e^{-it \mathfrak{B}} v_{0}\right\|_{H^1} &\leq C	\left\|\mathfrak{M}^{1+2s_1}|\textmd{Re}(v)|^{\alpha-1}\textmd{Re}(v)\right\|_{L^{p^\prime}([0,T];W^{{1},{q^\prime}})}\\
	&\leq CM^\alpha \leq M,
\end{aligned}	\\		
	\begin{aligned} \label{loc2}
			\sup _{\left[0, T\right]}\left\|\mathscr{T}	(v)-\mathscr{T}	(w)\right\|_{L^2(\mathbb{R}^{d})}&+\left\|\mathscr{T}	(v)-\mathscr{T}	(w)\right\|_{L^{p}([0,T];L^{q}(\mathbb{R}^{d}))}\\
		&\leq C
		\left\|\left\|v-w\right\|_{L^{q}}	\left(\left\|v\right\|^{\alpha-1}_{W^{{1},{q}}}+\left\|w\right\|^{\alpha-1}_{W^{{1},{q}}}\right)\right\|_{L^{p^\prime}[0,T]}\\
		&\leq C 
		M^{\alpha-1}\left\|v-w\right\|_{L^{p}([0,T];L^{q})}\leq \frac 12\left\|v-w\right\|_{L^{p}([0,T];L^{q})}.
	\end{aligned}	
	\end{align}
	As a consequence, we obtain that
	\begin{align*}
		\mathscr{T} (X_T)\subset X_T, \quad 	\mathfrak{d}\left(	\mathscr{T}(v), 	\mathscr{T}(w)\right)\leq \frac12 	\mathfrak{d}(v, w) 
	\end{align*}
and  $ \mathscr{T}(v) $ is a contraction map on $ (X_T,\mathfrak{d}) $.
Uniqueness   follows a similar process.

Next, we prove the continuity by utilizing the method in \cite{2003semilinear}. 
	Since $ \|v_{0,n}-v_0\|_{H^{1}}\rightarrow 0 $ and $ \|v_{0,n}\|_{H^{1}}\leq 2\|v_0\|_{H^{1}} $ for $ n\rightarrow \infty $, we can utilize Proposition \ref{propini} and the analysis in \eqref{loc} to find  $ N_1 $ and small $ T^{\prime}>0 $ such that, when $ n\geq N_1 $, one has
	\begin{align*}
			\left(\int_{0}^{T^{\prime}}\left\|e^{-it \mathfrak{B}}v_{0}\right\|_{q}^{p} d t\right)^{1 / p}<\varepsilon, 	~	\left(\int_{0}^{T^{\prime}}\left\|e^{-it \mathfrak{B}}v_{0,n}\right\|_{q}^{p} d t\right)^{1 / p}<\varepsilon\\
			\left\|\mathscr{T}	(v)\right\|_{L^{p}([0,T^{\prime}];W^{1,q})} \leq M,~\left\|\mathscr{T}	(v_n)\right\|_{L^{p}([0,T^{\prime}];W^{1,q})} \leq M,
	\end{align*}
	where $ M $ is the one set as in \eqref{loc0}. Noticing that
	\begin{align} \label{loc3}
		v_n(t)-v(t)= e^{-it \mathfrak{B}} (v_{0,n}-v_0) -i\beta\int_0^t 	e^{-i(t-\tau) \mathfrak{B}} \mathfrak{M}(|\textmd{Re}(v_n)|^{\alpha-1}\textmd{Re}(v_n)-|\textmd{Re}(v)|^{\alpha-1}\textmd{Re}(v))d \tau,
	\end{align}
 one can  discuss similarly as \eqref{loc2} to show
	\begin{align*}
		\sup _{\left[0, T^{\prime}\right]}\left\|	v_n-v	\right\|_{L^2}+\left\|	v_n-v\right\|_{L^{p}([0,T^{\prime}];L^{q})}& \leq C \left\|	v_{0,n}-v_0	\right\|_{L^2(\mathbb{R}^{d})}+C 
		M^{\alpha-1}\left\|v-w\right\|_{L^{p}([0,T^{\prime}];L^{q})}
	\end{align*} 
	and thus 
	\begin{align} \label{loc00}
		\sup _{\left[0, T^{\prime}\right]}\left\|	v_n-v	\right\|_{L^2}+ \frac 12\left\|	v_n-v\right\|_{L^{p}([0,T^{\prime}];L^{q})}& \leq C \left\|	v_{0,n}-v_0	\right\|_{L^2(\mathbb{R}^{d})}.
	\end{align}
Noticing that $ \nabla $ can commute with $ e^{-it \mathfrak{B}} $ and $ \mathfrak{M}  $, we can let  $ \nabla $ act on \eqref{loc3} to get
	\begin{align*} 
	\nabla(v_n-v)= &e^{-it \mathfrak{B}}\nabla(v_{0,n}-v_0) -i\beta\int_0^t 	e^{-i(t-\tau) \mathfrak{B}} \mathfrak{M}\alpha\left(|\textmd{Re}(v_n)|^{\alpha-1}(\nabla\textmd{Re}(v_n)-\nabla\textmd{Re}(v))\right)d \tau\\
&-i\beta\int_0^t 	e^{-i(t-\tau) \mathfrak{B}} \mathfrak{M}\alpha\left((|\textmd{Re}(v_n)|^{\alpha-1}-|\textmd{Re}(v)|^{\alpha-1})\nabla\textmd{Re}(v)\right)d \tau.
	\end{align*}
Then, similarly discussion as \eqref{loc2} shows
\begin{align} \label{loc4}
	\sup _{\left[0, T^{\prime}\right]}\left\|\nabla(v_n-v)	\right\|_{L^2}&+\left\|		\nabla(v_n-v)\right\|_{L^{p}([0,T^{\prime}];L^{q})}\nonumber\\
	& \leq C \left\| \nabla(v_{0,n}-v_0)	\right\|_{L^2(\mathbb{R}^{d})}+C 
	M^{\alpha-1}\left\|		\nabla(v_n-v)\right\|_{L^{p}([0,T^{\prime}];L^{q})}\\
	&+C\left\|\left(|\textmd{Re}(v_n)|^{\alpha-1}-|\textmd{Re}(v)|^{\alpha-1}\right)\nabla\textmd{Re}(v)\right\|_{L^{p^\prime}([0,T^{\prime}];L^{q^\prime})}. \nonumber
\end{align} 
Noticing $ 	C 
M^{\alpha-1} \leq \frac12 $, we can obtain the continuity by combining \eqref{loc00} and \eqref{loc4}, if one can prove
\begin{align}  \label{loc5}
\lim\limits_{n\rightarrow \infty}	\left\|\left(|\textmd{Re}(v_n)|^{\alpha-1}-|\textmd{Re}(v)|^{\alpha-1}\right)\nabla\textmd{Re}(v)\right\|_{L^{p^\prime}([0,T^{\prime}];L^{q^\prime})}=0.
\end{align}
Actually, if \eqref{loc5} fails, there exist a $ \epsilon_1>0 $ and  a subsequence $ \{v_{n_j}\}_{j\geq 1} $ such that
   \begin{align}  \label{loc6}
   		\left\|\left(|\textmd{Re}(v_{n_j})|^{\alpha-1}-|\textmd{Re}(v)|^{\alpha-1}\right)\nabla\textmd{Re}(v)\right\|_{L^{p^\prime}([0,T^{\prime}];L^{q^\prime})}\geq \epsilon_1  \text{ for } j\geq 1.
   \end{align}
However, we can use \eqref{loc00} to extract a subsequence, which is still denoted by $ \{v_{n_j}\}_{j\geq 1} $, such that as $ j\rightarrow \infty $, $  v_{n_j}\rightarrow v \text{ a.e. on }   [0,T^\prime]\times \mathbb{R}^d$ and   $  \textmd{Re}(v_{n_j})\rightarrow \textmd{Re}(v) \text{ a.e. on }   [0,T^\prime]\times \mathbb{R}^d$. Since $ \left(|\textmd{Re}(v_{n_j})|^{\alpha-1}-|\textmd{Re}(v)|^{\alpha-1}\right)\nabla\textmd{Re}(v) \in L^{p^\prime}([0,T^{\prime}];L^{q^\prime}) $, we obtain from the dominated convergence a contradiction with \eqref{loc6}.

	\begin{flushright}
		$\square$
	\end{flushright}

	\noindent{\bf{Proof of Theorem  \ref{th1}:}}
	
	Recall 
	\begin{align}
		u=&\frac{1}{2}\left(v+\bar{v}\right),~~
		(-\Delta)^{-\frac12}(1-\Delta)^{-\frac12}\partial_{t} u =\frac{i}{2}\left(\bar{v}-v\right) ,
	\end{align}
	thus, the transformation
	$$
	\begin{aligned}
		(v, \overline{v}) & \longrightarrow  (u,	(-\Delta)^{-\frac12}u_t) \\
		C\left([0, T]; H^{s}\right) \times C\left([0, T]; H^{s}\right) & \longrightarrow C\left([0, T]; H^{s}\right) \times C\left([0, T]; H^{s-1}\right)
	\end{aligned}
	$$
	is Lipschitz continuous and  Lemma \ref{le2}-\ref{le6} can yield the results.
	\begin{flushright}
		$\square$
	\end{flushright}

	\subsection{Global wellposedness and finite-time blowup} \label{s2.2}~
	
	As mentioned in the Introduction, despite  only giving the conclusion in dimension one, Liu \cite{liu1995instability} proved the preliminary Lemmas for all dimensions and the proof of global wellposedness, finite blowup and instability doesn't depend on the dimension. Thus, the conclusion can be applied to higher dimensions. Here, we give a proof based on the analysis of Wang et al. \cite{xmail}, Liu \cite{liu1995instability}. Compared with those of Liu \cite{liu1995instability}, our proof is more concise and the results of conclusion $ (i), (ii) $ are  better. We firstly present the preliminary Lemmas.
	
	\begin{lem} \label{le7} \cite{xmail}
		Assume that $y=y(t)$ is a continuous function. The constants $s, C_{1}$ $C_{2}$ satisfy $s>1, C_{2}>0,$ and $0<C_{1}<\frac{s-1}{s}\left(\frac{1}{C_{2} s}\right)^{\frac{1}{s-1}}$ are constants. If $0 \leq y(t) \leq$
		$C_{1}+C_{2} y(t)^{s},$ for any $t \geq 0,$ then there exist constants $y_{1}, y_{2},$ satisfying $0<y_{1}<$ $\left(\frac{1}{C_{2} s}\right)^{\frac{1}{s-1}}<y_{2}<\infty,$ such that
		\begin{equation}
			0\leq y(t) \leq y_{1}, \quad \text { if } \quad y(0)<y_{0}=\left(\frac{1}{C_{2} s}\right)^{\frac{1}{s-1}} \label{glob3}
		\end{equation}
		or
		\begin{equation}
			y_{2} \leq y(t)<\infty, \quad \text { if } \quad y(0)>y_{0}=\left(\frac{1}{C_{2} s}\right)^{\frac{1}{s-1}} \label{glob4}
		\end{equation}
		for any $t \geq 0$.
	\end{lem}
	
	\begin{lem} \label{le8} \cite{levine1974instability,xmail}
		Suppose that for $t \geq t_{0}, t_{0}$ is a positive constant, a nonnegative, twice-differentiable function $\phi(t)$ satisfies the inequality
		\begin{equation}
			\phi^{\prime \prime} \phi-(\eta+1)\left(\phi^{\prime}\right)^{2} \geq 0 \label{blow1}
		\end{equation}
		where $\eta>0$ is a constant. If $\phi\left(t_{0}\right)>0$ and $\phi^{\prime}\left(t_{0}\right)>0,$ then $\phi(t) \rightarrow \infty$ as $t \rightarrow t_{1} \leq$ $T_{0}=\phi\left(t_{0}\right) /\left(\eta \phi^{\prime}\left(t_{0}\right)\right)+t_{0}$.
	\end{lem}
	\begin{prop} \label{propg} Let $ \eta_d, E(u),  \varphi_d, C_\ast $ be the notation defined in \eqref{nota}, \eqref{ground},  \eqref{nota1}, $ \alpha $ satisfies \eqref{ch4_alpha0}, then
		\begin{align}
			&\eta_d= E(\varphi_d), \label{ass4} \\
			C_\ast= &C_{\alpha,d}= \|\varphi_d\|_{H^1(\mathbb{R}^d)}^{-\frac{\alpha-1}{\alpha+1}}.  \label{ass5} 
		\end{align}
	\end{prop}
	Proposition \ref{propg} is the combination of Theorem 2.6 and Corollary 2.5 in 	Liu \cite{liu1995instability}.

	\noindent{\bf{Proof of Theorem  \ref{th2}:}}
	Obviously, due to \eqref{ch4_est1}, the solution $ (u,u_t) \in H^1(\mathbb{R}^d)\times
	\dot{H}^{-1}(\mathbb{R}^d)) $ in Theorem \ref{th2} is global if $ \beta=1 $. When $ \beta=-1, $ we 
	let \[
	y(t)=\|u\|_{H^{1}}^{2},
	\]
	then \eqref{ch4_est1} shows
	\[
	\begin{aligned}
		y(t) & \leq 2 \mathcal{E}(0)+\frac{2}{\alpha+1}\|u\|_{L^{\alpha+1}}^{\alpha+1} \\
		& \leq 2\mathcal{E}(0)+\frac{2}{\alpha+1}C_{*}^{\alpha+1} y(t)^{\frac{\alpha+1}{2}}.
	\end{aligned}
	\]
	where $ C_\ast $ is the  best constant for $ \| u\|_{L^{\alpha+1}} \leq  C_\ast	\| u\|_{H^1} $.
	Combining with \eqref{assum1} $ ~\mathcal{E}(0)<\frac{\alpha-1}{2(\alpha+1)}C_\ast^{-\frac{2(\alpha+1)}{\alpha-1}},$  $ \|u_0\|_{H^1}<C_\ast^{-\frac{\alpha+1}{\alpha-1}}$
	and Lemma \ref{le7}, one has
	\begin{equation}
		\|u(t)\|_{H^1}^2<C_\ast^{-\frac{2(\alpha+1)}{\alpha-1}} \label{glob5},
	\end{equation}
	for $t>0$.
	It follows that
	\begin{align}
		&\frac{1}{\alpha+1}\|u\|_{L^{\alpha+1}}^{\alpha+1} \leq \frac{1}{\alpha+1}C_\ast^{\alpha+1}\|u\|_{H^1}^{\alpha+1}\nonumber\\
		&=\frac{1}{\alpha+1}C_\ast^{\alpha+1}\|u\|_{H^{1}}^{\alpha-1}\|u\|_{H^{1}}^{2}<\frac{1}{\alpha+1}\|u\|_{H^{1}}^{2},\label{glob6}
	\end{align}
	which combined with \eqref{ch4_est1} shows
	\begin{align}
		\frac{\alpha-1}{2(\alpha+1)}\|u\|_{H^{1}}^{2}
		+\frac 12\|(-\Delta)^{-\frac12}u_t\|_{L^2}^2 \leq   \mathcal{E}(0). \label{glob8}
	\end{align}
	This validates $ T_{max}=\infty $ under premise \eqref{assum1}.

	Next, we discuss the finite-time blowup in conclusion $ (ii) $ and
	we'll prove the argument by contradiction. Suppose that $T_{\max }=+\infty ,$ and for some $0<T<\infty$
	we can set
	\begin{equation}
		\phi(t)=\left\|(-\Delta)^{-\frac{1}{2}} u\right\|^{2}_{L^2}, \label{blow2}
	\end{equation}
	then,
	\begin{align}
		\phi^{\prime}(t) &=2\left((-\Delta)^{-\frac{1}{2}} u,(-\Delta)^{-\frac{1}{2}} u_{t}\right), \quad \forall t \in[0, T],
	\end{align}
	and
	\begin{align}
		\frac{1}{4}\left[\phi^{\prime}(t)\right]^{2} \leq \phi(t)\left\|(-\Delta)^{-\frac{1}{2}} u_{t}\right\|^{2}_{L^2}, \quad \forall t \in[0, T]. \label{blow6}
	\end{align}
	It follows the equation \eqref{ch4_1} that
	\begin{align}
		\phi^{\prime \prime}(t) &=2\left\langle(-\Delta)^{-1} u_{t t}, u\right\rangle_{H^{-1}, H_{0}^{1}}+2\left\|(-\Delta)^{-\frac{1}{2}} u_{t}\right\|^{2}_{L^2} \nonumber\\
		&=-2\|u\|_{H^{1}}^{2}+2 \int_{\mathbb{R}^{d}}  |u|^{\alpha+1} \mathrm{d} x+2\left\|(-\Delta)^{-\frac{1}{2}} u_{t}\right\|^{2}_{L^2}. \label{blow3}
	\end{align}
	Noticing that \eqref{ch4_est1} can be written into
	\begin{align*}
		2\int_{\mathbb{R}^d}  |u|^{\alpha+1} dx =-2(\alpha+1)\mathcal{E}(0)+(\alpha+1)\|u\|_{H^1}^2+(\alpha+1)\|(-\Delta)^{-\frac12}u_t\|_{L^2}^2,
	\end{align*}
	one has
	\begin{align}
		\phi^{\prime \prime}(t) &=(\alpha-1)\|u\|_{H^{1}}^{2}-2(\alpha+1)\mathcal
		{E}(0)+(\alpha+3)\left\|(-\Delta)^{-\frac{1}{2}} u_{t}\right\|^{2}_{L^2}. \label{blow7}
	\end{align}
	By Lemma \ref{le7} and premise \eqref{assum2}, we know there exists a constant $ y_2>0 $ such that
	\begin{equation}
		\|u(t)\|_{H^1}^2\geq y_2>C_\ast^{-\frac{2(\alpha+1)}{\alpha-1}},
	\end{equation}
	and
	\begin{align*}
		\phi^{\prime \prime}(t) &\geq (\alpha-1)(\|u\|_{H^{1}}^{2}-C_\ast^{-\frac{2(\alpha+1)}{\alpha-1}})\\
		&\geq  (\alpha-1)(y_2-C_\ast^{-\frac{2(\alpha+1)}{\alpha-1}})t>0,
	\end{align*}
	where 
	$$\mathcal{E}(0)<\frac{\alpha-1}{2(\alpha+1)}C_\ast^{-\frac{2(\alpha+1)}{\alpha-1}}.$$
	Therefore,
	\begin{align} \label{blow8}
		\phi^{\prime}(t) >\phi^{\prime}(0)+(\alpha-1)(y_2-C_\ast^{-\frac{2(\alpha+1)}{\alpha-1}})t.
	\end{align}
	For sufficiently large $t_0$ satisfying  $t_0< T  $, we have $\phi^{\prime}(t_0)>0$ and as $\phi(t_0)\geq 0$, we have $\phi(t)>0$ on $[t_0,T]$.
	
	On the other hand, \eqref{blow6} and \eqref{blow7} show
	\begin{align} \label{blow9}
		&0<\phi(t) \phi^{\prime \prime}(t)-\frac{1}{4}\left(\alpha+3\right) \phi^{\prime}(t)^{2}\nonumber\\
		&=\phi(t) \phi^{\prime \prime}(t)-\left(1+\frac{1}{4}(\alpha-1)\right) \phi^{\prime}(t)^{2}.
	\end{align}
	As a consequence, the Lemma \ref{le8} validates that there exists $t_{1} \leq T_{0}=\dfrac{\phi\left(t_{0}\right)}{\frac{1}{4}(\alpha-1) \phi^{\prime}\left(t_{0}\right)}+t_{0}$ such that
\begin{align} \label{goalblow}
		\lim _{t \to t_{1}} \phi(t)=+\infty.
\end{align}
	As 
	\begin{align*}
		\|(-\Delta)^{-\frac{1}{2}} u\|_{L^2}\leq\|(-\Delta)^{-\frac{1}{2}} u_{0}\|_{L^2}+\int_{0}^{t}\|(-\Delta)^{-\frac{1}{2}} u_{\tau} \|_{L^2}\mathrm{d} \tau, \quad \text { a.e. } t \in[0, T] .
	\end{align*} 
	we have 
	\begin{align*}
		\int_{0}^{t_1}\|(-\Delta)^{-\frac{1}{2}} u_{\tau} \|_{L^2}\mathrm{d} \tau=\infty
	\end{align*}
	and this implies there exists a sequence $ \{t_n\}, 0<t_n<t_1 $, such that 
	\begin{align*}
		\lim _{t_n \to t_{1}} \|(-\Delta)^{-\frac{1}{2}} u_{t} \|_{L^2}=\infty,
	\end{align*} 
	which contradicts to the assumption $T_{max}=\infty$.  Thus, the  solution $u$ will blowup under premise \eqref{assum2} and 
	\begin{align*}
		\lim _{t \to T_{max}^-} (\|u \|_{H^1}+\|(-\Delta)^{-\frac{1}{2}} u_{t} \|_{L^2})=\infty,
	\end{align*}
	which combined with \eqref{ch4_est1} gives
	\begin{align*}
		\lim _{t \to T_{max}^-} \|u \|_{L^{p+1}}=	\lim _{t \to T_{max}^-} \|u \|_{H^1}=\infty.
	\end{align*}

	Conclusion $ (iii) $ follows Theorem 4.3 and Lemma 4.4 of Liu \cite{liu1995instability}, as the proof method doesn't depend on the dimension. Thus, we omit the proof.
	\begin{flushright}
		$\square$
	\end{flushright}
	
	\begin{rem}
		The same  conclusion of $ (i), (ii) $  in Theorem \ref{th2} can be obtained if \eqref{assum1}, \eqref{assum2} are replaced with $   \beta=-1,\mathcal{E}(0)<\eta_d , R(u_0) >0 $
		and $ \mathcal{E}(0)<\eta_d , R(u_0) <0 $ respectively.	These correspond to Theorem 3.4, 4.2 of Liu \cite{liu1995instability}. We remark that our results are better, as we have
		\begin{align}
			&\eta_d = \frac{\alpha-1}{2(\alpha+1)}C_\ast^{-\frac{2(\alpha+1)}{\alpha-1}}. \label{ass1}\\
			\text{ If  }  \mathcal{E}(0)<\eta_d , R(u_0) >0,  &\text{    then }  \mathcal{E}(0)<\frac{\alpha-1}{2(\alpha+1)}C_\ast^{-\frac{2(\alpha+1)}{\alpha-1}} \|u_0\|_{H^1}<C_\ast^{-\frac{\alpha+1}{\alpha-1}}. \label{ass2}\\
			\text{ If  } & R(u_0) <0,  \text{ then  }\|u_0\|_{H^1}>C_\ast^{-\frac{\alpha+1}{\alpha-1}} . \label{ass3} 
		\end{align}
		In fact, multiplying \eqref{ground} with $ \varphi_d $ and integrating on $ \mathbb{R}^d $ gives $ \|\varphi_d\|_{H^1}^2=\|\varphi_d\|_{L^{\alpha+1}}^{\alpha+1} $. Then, \eqref{ass1} follows \eqref{ass4}, \eqref{ass5} that
		\begin{align*}
			\eta_d&=\frac12\|\varphi_d\|_{H^1}^2-
			\frac{1}{\alpha+1}\|\varphi_d\|_{L^{\alpha+1}}^{\alpha+1} \\
			&=\frac{\alpha-1}{2(\alpha+1)}\|\varphi_d\|_{H^1}^2=\frac{\alpha-1}{2(\alpha+1)}C_\ast^{-\frac{2(\alpha+1)}{\alpha-1}}
		\end{align*}
		and this  immediately shows the equivalence between $   \mathcal{E}(0)<\eta_d $ and $ \mathcal{E}(0)<\frac{\alpha-1}{2(\alpha+1)}C_\ast^{-\frac{2(\alpha+1)}{\alpha-1}} $. For \eqref{ass2},
		$ \mathcal{E}(0)<\frac{\alpha-1}{2(\alpha+1)}C_\ast^{-\frac{2(\alpha+1)}{\alpha-1}} $ and  $ R(u_0) >0  $, i.e. $ -\|u_0\|_{H^1}^2<-\|u_0\|_{L^{\alpha+1}}^{\alpha+1}  $ give 
		\begin{align*}
			\frac{\alpha-1}{2(\alpha+1)}C_\ast^{-\frac{2(\alpha+1)}{\alpha-1}} >\frac12\|u_0\|_{H^1}^2-
			\frac{1}{\alpha+1}\|u_0\|_{L^{\alpha+1}}^{\alpha+1}>\frac{\alpha-1}{2(\alpha+1)}\|u_0\|_{H^1}^2, 
		\end{align*}
		thus one has $ \|u_0\|_{H^1}^2<C_\ast^{-\frac{2(\alpha+1)}{\alpha-1}} $. \eqref{ass3} follows immediately by  $ R(u_0) <0  $ and the definition of $ C_{\ast} $
		\begin{align*}
			\|u_0\|_{H^1}^2<\|u_0\|_{L^{\alpha+1}}^{\alpha+1} \leq C_\ast^{\alpha+1} 	\|u_0\|_{H^1}^{\alpha+1},
		\end{align*}
		i.e.  $ 	\|u_0\|_{H^1}>C_\ast^{-\frac{\alpha+1}{\alpha-1}}  $.  
	\end{rem}

\begin{rem}
Consider the ``good" Boussinesq equation for all dimensions		
		\begin{align} \label{gB}
				\left\{\begin{array}{l}
					\partial_{t}^{2} u-\Delta u+\Delta^2 u+\Delta(u^2)=0, \quad (t,x)\in\mathbb{R}\times \mathbb{R}^d\\
				u(0, x)=u_{0}(x), ~~u_{t}(0, x)=u_{1}(x), 
			\end{array}\right.	
	\end{align}
which has conservation law 
\begin{align} \label{gBcons}
		\mathcal{E}_1(t): &=\frac 12\int_{\mathbb{R}^d}\left( |(-\Delta)^{-\frac12} u_t|^2 +u^2+|\nabla u|^2\right)-\frac13\int_{\mathbb{R}^d}u^3=  \mathcal{E}_1(0).
\end{align}
 Though the sign for term $ \frac13\int_{\mathbb{R}^d} u^3 $ in \eqref{gBcons} is undetermined,  the ``good" Boussinesq equation \eqref{gB} is somewhat similar with the focusing generalized Boussinesq equation \eqref{ch4_1} with $ \alpha=2 $ and it is easy to verify that Theorem \ref{th1} and the proof in Subsection \ref{s2.2} work for \eqref{gB} as well. For instance, if we can obtain 
 	\begin{align*}
	\lim _{t \to T_{max}^-} (\|u \|_{H^1}+\|(-\Delta)^{-\frac{1}{2}} u_{t} \|_{L^2})=\infty,
\end{align*}
then \eqref{gBcons}  shows 
\begin{align*}
\lim _{t \to T_{max}^-}\frac13\int_{\mathbb{R}^d}u^3	=\lim _{t \to T_{max}^-}\left[\frac 12\int_{\mathbb{R}^d}\left( |(-\Delta)^{-\frac12} u_t|^2 +u^2+|\nabla u|^2\right)- \mathcal{E}_1(0)\right]= \infty,
\end{align*}
which combined with
$ \int_{\mathbb{R}^d}u^3 \leq\int_{\mathbb{R}^d}|u|^3 \leq  C_{2,d}^3	\|u_0\|_{H^1}^3 $ gives
\begin{align*}
	\lim _{t \to T_{max}^-} \|u \|_{L^{3}}=	\lim _{t \to T_{max}^-} \|u \|_{H^1}=\infty.
\end{align*}
So we have following theorem.
\begin{thm}~\\
	$ (i) $ 		The initial value problem \eqref{gB} with $(u, (-\Delta)^{-\frac12}u_t)$ in $H^{s} (\mathbb{R}^{d}) \times H^{s-1} (\mathbb{R}^{d})$ is locally well-posed  for
	\begin{align*}  
		\begin{array}{ll}
			s=0, &  \text { if } d\leq 4, \\
			s = 1 , & \text { if } 
		 d\leq 6.
		\end{array}
	\end{align*}
	$ (ii) $ 		Suppose $ (u_0,u_1)\in H^{1}(\mathbb{R}^{d})\times \dot{H}^{-1}(\mathbb{R}^{d}),d\leq5 $, then 	the initial value problem \eqref{gB} with $(u, u_t)$  in $ H^{1} (\mathbb{R}^{d})\times \dot{H}^{-1} (\mathbb{R}^{d})$ is globally well-posed if
	\begin{align*}  
\mathcal{E}_1(0)<\frac{1}{6C_{2,d}^6},\quad  \|u_0\|_{H^1}<C_{2,d}^{-3}, 	
	\end{align*}
and blows up in finite time if 
\begin{align*}
u_0\in \dot{H}^{-1}(\mathbb{R}^{d}),~	\mathcal{E}(0)<\frac{1}{6C_{2,d}^6},~~ \|u_0\|_{H^1}>C_{2,d}^{-3}.	
\end{align*} 	
\end{thm}

\end{rem}
	
	\section{Small initial data scattering} \label{s3.}~
	
	Similar with the arguments of local wellposedness, we can follow corresponding theory of \eqref{nsl1} to arrive at the results of small initial data scattering.
	
	\noindent{\bf{Proof of Theorem  \ref{th3}:}}
	
	Define
	\begin{align}
		\|v\|_{S^{s}\left(I \times R^{d}\right)}:=\underset{(p, q) \text { admissible }}{\sup }\|\left\langle\nabla\right\rangle^s v\|_{L_{t}^{p} L_{x}^{q}\left(I \times R^{d}\right)} \label{S1n}
	\end{align}	
	and let  $ q_1^\prime= (d+2)/2,$ $ 
	p_1^\prime= 2(d+2)/(d+6) $, $ s_c=d / 2-2 / (\alpha-1)$. Noticing that 
	the assumption  gives
	\begin{align}
		\|v_0\|_{H^1}&\leq\|u_0\|_{H^1}+\| \left(-\Delta\right)^{-\frac12}\left(1-\Delta\right)^{-\frac12}  u_1\|_{H^1} \nonumber\\
		&=\|u_0\|_{H^1}+\|  u_1\|_{\dot{H}^{-1}}
		\leq  \varepsilon, \label{ch4_ini}
	\end{align}
	where $ \varepsilon $ is a small enough constant,   then by \eqref{ch4_inte2}, \eqref{stype4}, \eqref{stype5}, \eqref{nonlinear2} and Sobolev embedding inequality, one has
	\begin{equation} \label{scat1}
		\begin{aligned}
			\left\|\mathscr{T} v\right\|_{S^{1}\left((0,t) \times R^{d}\right)}&	\leq C\left\|v_0\right\|_{H^{1}}+C \left\||v|^{\alpha}\right\|_{L^{q_{1}^{\prime}}([0,t];W^{1,p_{1}^\prime})  } \\
			& \leq C \varepsilon+C \left\|v\right\|_{L^\infty_tH^{1}_x }  \left\||v|^{\alpha-1}\right\|_{L^{\frac{d+2}{2}}([0,t];L^{\frac{d+2}{2}})  }\\
			& \leq C \varepsilon+C \left\|v\right\|_{L^\infty_tH^{1}_x } \left\|v\right\|^{\alpha-1}_{L^{\frac{(\alpha-1)(d+2)}{2}}\left([0,t];W^{1,l_1}\right)  }\\
			& \leq C \varepsilon+C	\|v\|_{S^{1}\left((0,t) \times R^{d}\right)}^\alpha,
		\end{aligned}
	\end{equation}
	where $ l_1=\frac{2(\alpha -1)d(d+2)}{(\alpha-1) d(d+2)-8} $, 
	$ ((\alpha-1)(d+2)/2,l_1) $ is admissible,  
	$ \alpha $ satisfies \eqref{ch4_alpha1}, $ 0\leq s_c \leq 1 $ and we used Sobolev embedding relationship 
	$ W^{1,l_1}(\mathbb{R}^d)\hookrightarrow W^{s_c,l_1}(\mathbb{R}^d)\hookrightarrow L^{\frac{(\alpha-1)(d+2)}{2}}(\mathbb{R}^d) $ 
	(see \cite{di2012hitchhikers}).
	Then, standard bootstrap argument shows
	\begin{align}
		\|v\|_{S^{1}\left(\mathbb{R} \times R^{d}\right)} \leq C. \label{scat1'}
	\end{align}
	Let
	\begin{align} \label{v+}
		v^{+}=v_{0}-i \beta\int_{0}^{+\infty} e^{i t^{\prime} \mathfrak{B}}\mathfrak{M}|\textmd{Re}v|^{\alpha-1}\textmd{Re} vd t^{\prime},
	\end{align}
	it follows \eqref{ch4_inte1} that
	\begin{align}
		v(t)-e^{-it \mathfrak{B}}v^{+}=i\beta \int_{t}^{+\infty} e^{-i\left(t-t^{\prime}\right) \mathfrak{B}}\mathfrak{M}|\textmd{Re}v|^{\alpha-1}\textmd{Re}v d t^{\prime}  \label{ch4_scat2}
	\end{align}
	and  Strichartz estimate \eqref{stype4},  \eqref{stype5}, process of \eqref{scat1}  that
	\begin{equation}
		\begin{aligned}
			\left\|	v(t)-e^{-it \mathfrak{B}}v^{+}\right\|_{H^{1}} & \lesssim  \left\||v|^{\alpha}\right\|_{L^{q_{1}^{\prime}}(t,\infty;W^{1,p_{1}^\prime})  } \\
			& \lesssim \left\|v\right\|_{L^\infty_tH^{1}_x } \left\|v\right\|^{\alpha-1}_{L^{\frac{(\alpha-1)(d+2)}{2}}\left(t,\infty;W^{s_c,l_1}\right)  }\\
			& \lesssim  \|v\|_{S^{1}\left((t,\infty) \times R^{d}\right)}^\alpha \leq C,
		\end{aligned}
	\end{equation}
	and thus
	\begin{align}
		\lim_{t \to \infty}	\left\|	v(t)-e^{-it \mathfrak{B}}v^{+}\right\|_{H^{1}}=0. \label{ch4_scat3}
	\end{align}
	Similarly,
	\begin{align}
		\lim\limits_{t\rightarrow - \infty}\left\|	v(t)-e^{-it \mathfrak{B}}v^{-}\right\|_{H^{1}}=0 \label{ch4_scat4}
	\end{align} 
	for
	\begin{align} \label{v-}
		v^{-}=v_{0}-i \beta\int^{0}_{-\infty} e^{i t^{\prime} \mathfrak{B}}\mathfrak{M}|\textmd{Re}v|^{\alpha-1}\textmd{Re} vd t^{\prime}.
	\end{align}
	Simple calculation on \eqref{v+}, \eqref{v-}, \eqref{ch4_s1} gives
	\begin{align*}
		\textmd{Re}(e^{-it \mathfrak{B}}v^{\pm})&= \cos(t \mathfrak{B})\textmd{Re} ( v_{0})+\sin(t \mathfrak{B})\textmd{Im} ( v_{0})+\beta\int_{I^{\pm}}  \sin((t-t^{\prime}) \mathfrak{B})\mathfrak{M}|\textmd{Re}v|^{\alpha-1}\textmd{Re} vd t^{\prime}\\
		&=\cos(t \mathfrak{B})u_0+\sin(t \mathfrak{B})\mathfrak{B}^{-1} u_1+\beta\int_{I^{\pm}}  \sin((t-t^{\prime}) \mathfrak{B})\mathfrak{M}|\textmd{Re}v|^{\alpha-1}\textmd{Re} vd t^{\prime},\\
		\textmd{Im}(e^{-it \mathfrak{B}}v^{\pm})&=- \sin(t \mathfrak{B})u_0+\cos(t \mathfrak{B})\mathfrak{B}^{-1} u_1-\beta\int_{I^{\pm}}  \cos((t-t^{\prime}) \mathfrak{B})\mathfrak{M}|\textmd{Re}v|^{\alpha-1}\textmd{Re} vd t^{\prime}, 
	\end{align*}
	where  $ I^+=[0,\infty], I^-=[-\infty,0] $ and 
	\begin{align} \label{cosB}
		\cos\left(t\mathfrak{B}\right)= \mathscr{F}_\xi^{-1} \cos\left( t|\xi|\sqrt{|\xi|^2+1} \right)\mathscr{F}_x,~ \sin\left(t\mathfrak{B}\right)=  \mathscr{F}_\xi^{-1} \sin\left( t|\xi|\sqrt{|\xi|^2+1} \right)\mathscr{F}_x.
	\end{align}
	Define $u_0^{\pm}:=\textmd{Re}(e^{-it \mathfrak{B}}v^{\pm}),~u_1^{\pm}:=\mathfrak{B}\textmd{Im}(e^{-it \mathfrak{B}}v^{\pm})$, then  \eqref{ch4_scat3}, \eqref{ch4_scat4} show \eqref{sc1}:
	\begin{align} \label{scat6}
		\begin{aligned}
			&	\lim _{t \rightarrow \pm \infty}\left\|\left(u(t), u_{t}(t)\right)-B(t)\left(\textmd{Re}(v^{\pm}), \textmd{Im}(v^{\pm})\right)\right\|_{\mathcal{H}^{1}}\\
			=&\lim _{t \rightarrow \pm \infty}\left\|\left(u(t), u_{t}(t)\right)-\left(u_0^{\pm}, u_1^{\pm}\right)\right\|_{\mathcal{H}^{1}}\\
			=& \lim _{t \rightarrow \pm \infty}	\left\|	v(t)-e^{-it \mathfrak{B}}v^{\pm}\right\|_{H^{1}}=0,
		\end{aligned}
	\end{align}
	where $ \mathcal{H}^{1}=H^1\times\dot{H}^{-1}  $, 
	$ 	B(t)= \left(\begin{array}{cc}
		\cos(t \mathfrak{B}) &  \sin(t \mathfrak{B})\\
		-\sin(t \mathfrak{B}) & \cos(t \mathfrak{B})
	\end{array}\right) $
	and $(u(t), \mathfrak{B}^{-1}u_t(t))= B(t)(\tilde{u}_0, \mathfrak{B}^{-1}\tilde{u}_1) $ solves the linear part of gBQ, i.e. 
	\begin{align*} 
		\left\{\begin{array}{l}
			\partial_{t}^{2} u-\Delta u+\Delta^2 u=0,\\
			u(0, x)=\tilde{u}_{0}(x), ~~u_{t}(0, x)=\tilde{u}_{1}(x).  
		\end{array}\right.
	\end{align*}
	\begin{flushright}
		$\square$
	\end{flushright}
	\begin{rem}
		The mass critical case $ \alpha=1+4/d $ in Theorem \ref{th3} can be improved as 
		\begin{align*}
			\lim _{t \rightarrow \pm \infty}\left\|\left(u(t), u_{t}(t)\right)-\left(u_0^{\pm}, u_1^{\pm}\right)\right\|_{\mathcal{L}^{2}}=0,
		\end{align*}
		for given $ (u_0, (-\Delta)^{-\frac12}u_1)\in L^{2}(\mathbb{R}^{d})\times H^{-1}(\mathbb{R}^{d}) $ small enough. Here,  $ \left\|\left(w_1,w_2\right)\right\|_{\mathcal{L}^{2}}:=(\|w_1\|_{L^{2}}^{2}+\|(-\Delta)^{-\frac12}w_2\|_{H^{-1}}^{2})^{\frac{1}{2}} $, $ \left(u_0^{\pm}, (-\Delta)^{-\frac12}u_1^{\pm}\right) \in L^{2}(\mathbb{R}^{d})\times H^{-1}(\mathbb{R}^{d})$  is the linear solution of \eqref{ch4_1}. In fact, we can follow the arguments of Theorem 1.17 in Dodson \cite{dodson2019defocusing} to obtain that \eqref{ch4_s1} scatters in $ L^2(\mathbb{R}^d) $  if $\| v_0\|_{  L^2(\mathbb{R}^d)}$ small enough. Then, similar deduction with \eqref{scat6} shows the conclusion.
		Likewise, following Theorem 3.1 of  \cite{dodson2019defocusing}, the energy critical case $ \alpha=1+4/(d-2),d\geq 3 $ can be improved as
		\begin{align*}
			\lim _{t \rightarrow \pm \infty}\left\|\left(u(t), u_{t}(t)\right)-\left(u_0^{\pm}, u_1^{\pm}\right)\right\|_{\mathcal{\dot{H}}^{1}}=0,
		\end{align*}
		for given $ (u_0, u_1)\in \dot{H}^{1}(\mathbb{R}^{d})\times H^{-1}(\mathbb{R}^{d}) $ small enough, where 
		$  \left\|\left(w_1,w_2\right)\right\|_{\mathcal{\dot{H}}^{1}}:=(\|w_1\|_{\dot{H}^{1}}^{2}+\|w_2\|_{H^{-1}}^{2})^{\frac{1}{2}} $ and $ \left(u_0^{\pm}, u_1^{\pm}\right) \in \dot{H}^{1}(\mathbb{R}^{d})\times H^{-1}(\mathbb{R}^{d})$  is the linear solution of \eqref{ch4_1}. 

	\end{rem}
	
	\section {Large initial data scattering for  radial defocusing case} ~\label{s4.}
	
	To show larger initial data scattering, it is classical to establish Morawetz-virial type estimate. Similar to the wave and Klein-Gordon equation, we define
	\begin{align*}	
		\qquad M_{a}(t):=-\int_{\mathbb{R}^{d}}\left((-\Delta)^{-\frac{1}{2}} u_{t}\right)\left(\nabla a \cdot \nabla(-\Delta)^{-\frac{1}{2}} u+\frac{1}{2} \Delta a(-\Delta)^{-\frac{1}{2}} u\right) d x,
	\end{align*}
	where $ a $ will be clarified later. See also \cite{nakanishi2001remarks} for more general form of the definition of such quality. Combining the bootstrap argument and Morawetz-virial type estimate, Dodson-Murphy \cite{dodson2017new} give a new proof of the scattering
	result in \cite{holmer2008sharp}  about the 
	cubic Schr\"odinger equation in dimension three, which can avoid the use of concentration compactness argument. 
	See also \cite{guo2020scattering,chen2021scattering} for similar argument for other dispersive models. The new  ingredient here is that we need to deal carefully with the operator $ (-\Delta)^{-\frac12}$. Explicitly, we need Lemma \ref{le4.1} to estimate the commutator between $ (-\Delta)^{\frac12}  $ and some good enough functions.
	
	\begin{prop}  Space-time estimate\label{prop4.1}
		
		Suppose $ 4/d+1\leq\alpha\leq(d+2)/(d-2) $, $ d\geq 3  $ and  the $ H^1(\mathbb{R}^d) $  solution $ u $ of \eqref{ch4_1} is global wellposed, then $ u $ satisfies
		\begin{align} \label{tsinte1}
			\int_{I} \int_{\mathbb{R}^d} |u|^{\alpha+1} dx dt\leq 	\left\{\begin{array}{l}
				 C|I|^{\frac{1}{1+\theta}}, \quad~~ \text{ if } d\geq 4, \\[5pt] 
		C\frac{|I|}{\log |I|},  \quad~~~ \text{ if } d=3 \text{ and } |I|\geq 2,
			\end{array}\right.
		\end{align}
		where $ \theta  $ is a positive constant less than 1 and  $ I $ is any time interval in $ \mathbb{R}$.
	\end{prop}
	
	\noindent{\bf{Proof of Proposition \ref{prop4.1}:}}
	
	 Denote  $w=(-\Delta)^{-\frac{1}{2}} u, w_{j}=\partial_{x_{j}} w, w_{k l}=\partial_{x_{k}} \partial_{x_{l}} w, a_{j k}=\partial_{x_{j}} \partial_{x_{k}} a$, then $ \partial_{x_{j}} w= - R_{j} u  $ and
	\begin{align} \label{Ma'}
		M_{a}&^{\prime}(t)= \int_{\mathbb{R}^{d}} a_{j k}\left(2 w_{j l} w_{k l}+w_{j} w_{k}\right) d x \nonumber\\
		&+\int_{\mathbb{R}^{d}} \frac{1}{4}\left(\Delta^{3} a-\Delta^{2} a\right) w^{2} -\frac{1}{2} \Delta^{2} a|\nabla w|^{2}-\frac{1}{2} \Delta a_{j k} w_{j} w_{k}+a_{j k l} w_{j} w_{k l} d x \\
		&+\int_{\mathbb{R}^{d}}|u|^{\alpha-1} u(-\Delta)^{\frac{1}{2}}\left(\nabla a \cdot \nabla(-\Delta)^{-\frac{1}{2}} u+\frac{1}{2} \Delta a(-\Delta)^{-\frac{1}{2}} u\right) d x. \nonumber
	\end{align}
	For  case $ d\geq 4 $, we choose $a(x)=R^{2} a_{0}(x / R),  R \geq 1$ to define the Morawetz quality $M_{R}$, where  $a_{0}(x)=\tilde{a}_{0}(|x|)$ and   $\tilde{a}_{0}(r)\in C^\infty$ is a smooth function in $\mathbb{R}$ such that
	\begin{equation} \label{a}
		\tilde{a}_{0}(r)=\left\{\begin{array}{ll}
			r^{2}, & r<1/2 \\
		 \gamma_1 r+\gamma_2 , & r>1
		\end{array}\right.
	\end{equation}
 $  \text{ for constants } \gamma_1>0, \gamma_2 \in \mathbb{R}  $,  and
	\begin{align}
		\tilde{a}_{0}^{\prime}(r) \geq 0, \quad \tilde{a}_{0}^{\prime \prime}(r) \geq 0, \quad\left(\partial_{r}^{2}+\frac{d-1}{r} \partial_{r}\right)^{2} \tilde{a} \leq 0, \quad \forall r>0. \label{a1}
	\end{align}
	Thus, one has
	\begin{align*}
		&\Delta^2 a=0  ~\text {    for     } ~ |x|<R/2,~~~~ \Delta^2 a\leq0  ~\text {    for     }  ~ R/2<|x|\leq R, 	\\
		& \left|\partial^{\alpha} a(x)\right| \lesssim_{\alpha} R|x|^{-|\alpha|+1} ~\text { for } ~|\alpha| \geq 1 ,~ R/2\leq|x|.
	\end{align*}	
Note that it is easy to prove  there is no nontrivial $ \tilde{a}_{0}(r) $ in $ d=3 $ that satisfies both \eqref{a} and \eqref{a1}. Since $\left(\left(a_{R}\right)_{j k}\right)$ is a nonnegative matrix, we have
	\begin{equation} \label{MR'}
		\begin{aligned}
			M_{R}^{\prime}(t) \geq & \int_{\mathbb{R}^{d}}|u|^{\alpha-1} u(-\Delta)^{\frac{1}{2}}\left(\nabla a_{R} \cdot \nabla(-\Delta)^{-\frac{1}{2}} u+\frac{1}{2} \Delta a_{R}(-\Delta)^{-\frac{1}{2}} u\right) d x \\
			&-\int_{\mathbb{R}^{d}} \frac{1}{4}\left|\Delta^{3} a_{R}\right| w^{2}+\frac{1}{2}\left(\left|\Delta^{2} a_{R}\right|+\sum_{j, k}\left|\Delta\left(a_{R}\right)_{j k}\right|\right)|\nabla w|^{2} d x \\
			&-\int_{\mathbb{R}^{d}} \sum_{j, k, l}\left|\left(a_{R}\right)_{j k l}\right| |w_{j}|| w_{k l} \mid d x  \\
			\geq & \int_{\mathbb{R}^{d}}|u|^{\alpha-1} u(-\Delta)^{\frac{1}{2}}\left(\nabla a_{R} \cdot \nabla(-\Delta)^{-\frac{1}{2}} u+\frac{1}{2} \Delta a_{R}(-\Delta)^{-\frac{1}{2}} u\right) d x  \\
			&-C \int_{|x| \geq \frac{R}{2}} \frac{R}{|x|^{5}}|w|^{2}+\frac{R}{|x|^{3}}|\nabla w|^{2}+\frac{R}{|x|^{2}}\left|\nabla w \| \nabla^{2} w\right| d x.
		\end{aligned}
	\end{equation}
	It follows  Hardy inequality that
	\begin{align*}
		\int_{|x| \geq \frac{R}{2}} \frac{R}{|x|^{5}}|w|^{2} \lesssim  \frac{1}{R^2} 	\int_{\mathbb{R}^d} |\nabla w|^{2} \lesssim \frac{1}{R^2} 	\|u\|^2_{L^2}. 
	\end{align*}
	Then, H\"older inequality and \eqref{MR'} give  
	\begin{align} \label{MR'1}
		M_{R}^{\prime}(t) \geq \int_{\mathbb{R}^{d}}|u|^{\alpha-1} u(-\Delta)^{\frac{1}{2}}\left(\nabla a_{R} \cdot \nabla(-\Delta)^{-\frac{1}{2}} u+\frac{1}{2} \Delta a_{R}(-\Delta)^{-\frac{1}{2}} u\right) d x 
		-\frac{C}{R}\|u\|_{H^{1}}^{2}.
	\end{align}
	For case $ d= 3 $, we  can set $ a(x) $ as those in $ d\geq 4 $ and $ 	\tilde{a}_{0} $ as
 \begin{equation} \label{a2}
	\tilde{a}_{0}(r)=\left\{\begin{array}{ll}
		r^{2}, & r<1/2 \\
		2 r, & r>1,
	\end{array}\right.
\end{equation}
also
$$
\tilde{a}_{1}^{\prime}(r) \geq 0, \quad \tilde{a}_{1}^{\prime \prime}(r) \geq 0, \quad \forall r>0.
$$
As $ \Delta^2 a=0  ~\text {    for     } ~ |x|<R/2  \text{ or   } |x|>R $,  the same analysis as those in $ d\geq 4 $
shows
	\begin{align} \label{MR'11}
\begin{aligned}
		M_{R}^{\prime}(t) \geq& \int_{\mathbb{R}^{3}}|u|^{\alpha-1} u(-\Delta)^{\frac{1}{2}}\left(\nabla a_{R} \cdot \nabla(-\Delta)^{-\frac{1}{2}} u+\frac{1}{2} \Delta a_{R}(-\Delta)^{-\frac{1}{2}} u\right) d x \\
	&-\frac{C}{R}\|u\|_{H^{1}(\mathbb{R}^{3})}^{2}-\frac{C}{R^2} \int_{R / 2 \leq|x| \leq R}|u|^{2} d x.
\end{aligned}
\end{align}

	If one can derive for  $  d\geq 3 $ and $ \text{for  some constant } c>0, \theta>0 $ that 
	\begin{align} \label{goal2}
		\begin{aligned}
			& \int_{\mathbb{R}^{d}}|u|^{\alpha-1} u(-\Delta)^{\frac{1}{2}}\left(\nabla a_{R} \cdot \nabla(-\Delta)^{-\frac{1}{2}} u+\frac{1}{2} \Delta a_{R}(-\Delta)^{-\frac{1}{2}} u\right) d x \\
			\geq & c \int_{\mathbb{R}^{d}}|u|^{\alpha+1} d x-\frac{C}{R^{\theta}}\left(\|u\|_{H^{1}}^{2}+\|u\|_{H^{1}}^{\alpha+1}\right) ,
		\end{aligned}
	\end{align}
we can utilize \eqref{MR'1} and \eqref{MR'11} to get \eqref{tsinte1}.  	Actually, for case $ d\geq 4 $,  \eqref{MR'1} and \eqref{goal2} show
	\begin{equation} \label{goal1}
		M_{R}^{\prime}(t) \geq c \int_{\mathbb{R}^{d}}|u|^{\alpha+1} d x-\frac{C}{R^\theta}\left(\|u\|_{H^{1}}^{2}+\|u\|_{H^{1}}^{\alpha+1}\right),~.
	\end{equation}
and the conservation law \eqref{ch4_est1}, the definition of $ a_R $ and Cauchy-Schwarz inequality give $ |M_{R}(t)| \lesssim R \|(-\Delta)^{\frac12}u_t\|_{L^2}\|u\|_{H^1}\lesssim  R\mathcal{E}(0)  $. Integrating \eqref{goal1} on any time interval $ I $, we have 
	\begin{align*}
		\int_I \int_{\mathbb{R}^{d}}|u|^{\alpha+1} d x dt\leq C R+\frac{CI}{R^\theta}.
	\end{align*}
	If $ |I|>1 $, taking $ R=|I|^{\frac{1}{1+\theta}} $, we obtain \eqref{tsinte1}.
	If $ |I|\leq 1 $, the \eqref{tsinte1} is trivial by the conservation law \eqref{ch4_est1}. For case $ d=3 $,  we can use \eqref{MR'11} and \eqref{goal2} to get
	\begin{align} \label{d3}
		\begin{aligned}
			M_{R}^{\prime}(t) &\geq c \int_{\mathbb{R}^{3}}|u|^{\alpha+1} d x-\frac{C}{R^\theta}\left(\|u\|_{H^{1}(\mathbb{R}^3)}^{2}+\|u\|_{H^{1}(\mathbb{R}^3)}^{\alpha+1}\right)-\frac{C}{R^2}\int_{R / 2 \leq|x| \leq R}|u|^{2} d x\\
			&\geq c \int_{\mathbb{R}^{3}}|u|^{\alpha+1} d x-\frac{C}{R^\theta}-\frac{C}{R^2} \int_{R / 2 \leq|x| \leq R}|u|^{2} d x.		
		\end{aligned}
	\end{align}
Integrating \eqref{d3} on any time interval $ I $ and similar discussion as before give
\begin{align} \label{d3_1}
	\int_{I} \int_{\mathbb{R}^{3}}|u|^{\alpha+1} d x d t \leq C\left(R+R^{-\theta}|I|\right)+C \int_{I} \int_{R / 2 \leq|x| \leq R} R^{-2}|u|^{2} d x d t,
\end{align}
where the $ C  $ in the right hand side is independent of $ R $.
Assume $ |I| \geq 2 $,   we can sum \eqref{d3_1} over $ R $ with $ R= 2^{N} \leq|I|^{1 /(1+\theta)}, N \in \mathbb{N}^\ast $
 and use Hardy inequality to get
\begin{align*}
	\sum_{R \in 2^{N}, R \leq|I|^{1 /(1+\theta)}} \int_{I} \int_{\mathbb{R}^{3}}|u|^{\alpha+1} d x d t&\leq 
	C|I|^{\frac{1}{1+\theta}}+C|I|+C\int_{I} \int_{\mathbb{R}^{3}} \frac{|u(x)|^{2}}{|x|^{2}} d x dt \\
		&\leq C|I|^{\frac{1}{1+\theta}}+C|I| 
\end{align*}
	and then obtain
	\begin{align*}
		\int_{I} \int_{\mathbb{R}^{3}}|u|^{\alpha+1} d x d t \leq C\frac{|I|}{\log |I|}, 
	\end{align*}
	which  is  \eqref{tsinte1}.
	
Now, we only need to prove \eqref{goal2}. 
	 Directly calculation shows that
	\begin{equation} \label{}
		\begin{aligned}
			& \int_{\mathbb{R}^{d}}|u|^{\alpha-1} u(-\Delta)^{\frac{1}{2}}\left(\nabla a_{R} \cdot \nabla(-\Delta)^{-\frac{1}{2}} u+\frac{1}{2} \Delta a_{R}(-\Delta)^{-\frac{1}{2}} u\right) d x \\
			=& \int_{\mathbb{R}^{d}}(-\Delta)^{-\frac{1}{2}}\left(|u|^{\alpha-1} u\right)(-\Delta)\left(\partial_{j} a_{R} \partial_{j}(-\Delta)^{-\frac{1}{2}} u+\frac{1}{2} \partial_{j j} a_{R}(-\Delta)^{-\frac{1}{2}} u\right) d x \\
			=&-\int_{\mathbb{R}^{d}}(-\Delta)^{\frac{1}{2}}\left(\partial_{j}(-\Delta)^{-\frac{1}{2}}|u|^{\alpha-1} u \partial_{j} a_{R}\right) u d x \\
			&-2 \int_{\mathbb{R}^{d}}(-\Delta)^{-\frac{1}{2}}\left(|u|^{\alpha-1} u\right) \partial_{j} \partial_{k} a_{R} \partial_{j} \partial_{k}(-\Delta)^{-\frac{1}{2}} u d x \\
			&-\int_{\mathbb{R}^{d}} \frac{1}{2}(-\Delta)^{-\frac{1}{2}}|u|^{\alpha-1} u \Delta a_{R}(-\Delta)^{\frac{1}{2}} u d x \\
			&+\int_{\mathbb{R}^{d}}(-\Delta)^{-\frac{1}{2}}|u|^{\alpha-1} u\left(-2 \Delta \partial_{j} a_{R} \partial_{j}(-\Delta)^{-\frac{1}{2}} u-\frac{1}{2} \Delta^{2} a_{R}(-\Delta)^{-\frac{1}{2}} u\right) \\
			:=& \textmd{I}+\textmd{II}+\textmd{III}+\textmd{IV} .
		\end{aligned}
	\end{equation}
	By H\"older inequality, Sobolev embedding inequality $ \left\|v\right\|_{L^{\frac{dp}{d-p}}(\mathbb{R}^d)} \lesssim \left\|\nabla v\right\|_{L^{p}(\mathbb{R}^d)} $ and interpolation inequality, we have
	\begin{equation} \label{IV}
		\begin{aligned}
			\left|\textmd{IV}\right|	&\lesssim \int_{|x| \geq \frac{R}{2}}\left|(-\Delta)^{-\frac{1}{2}}| u|^{\alpha-1} u \right| \left(\frac{R}{|x|^{2}}\left|\nabla(-\Delta)^{-\frac{1}{2}} u\right|+\frac{R}{|x|^{3}}\left|(-\Delta)^{-\frac{1}{2}} u\right|\right) d x \\
			&\lesssim \frac{1}{R}\left\|(-\Delta)^{-\frac{1}{2}}|u|^{\alpha-1} u\right\|_{L^{2}}\left(\left\|\nabla(-\Delta)^{-\frac{1}{2}} u\right\|_{L^{2}}+\left\|\frac{1}{|x|}(-\Delta)^{-\frac{1}{2}} u\right\|_{L^{2}}\right)\\
			&\lesssim \frac{1}{R}\left\||u|^{\alpha-1} u\right\|_{L^{\frac{2 d}{d+2}}}\|u\|_{L^{2}}\\
			&\lesssim \frac{1}{R}\left\|u\right\|^\alpha_{H^1}\|u\|_{L^{2}}\lesssim \frac{1}{R},
		\end{aligned}
	\end{equation}
	where we have used  $ \frac{2d\alpha}{2+d}  \leq\frac{2d}{d-2}$  or $\alpha\leq(d+2)/(d-2) $.
	
	To dispose the commutation between $ (-\Delta)^{\frac{1}{2}} $ and terms of $ a_{R}  $ in $ \textmd{I},\textmd{II},\textmd{III} $,  we need the following important lemma from \cite{calderon1965commutators}, where we give a short clarification
	in Appendix \ref{appb}.
	\begin{lem} \label{le4.1}
		Assume $\varphi \in C^{\infty}, \nabla \varphi \in L^{\infty}$, $1<p<\infty$. Then,
		\begin{align}
			\left\|[(-\Delta)^{\frac{1}{2}}, \varphi] f\right\|_{L^{p}( \mathbb{R}^d)  } \lesssim\|\nabla \varphi\|_{L^{\infty}( \mathbb{R}^d)} \left\|f\right\|_{L^{p}( \mathbb{R}^d)  }. \label{commu1}
		\end{align}	
	\end{lem}
	We also need the radial Sobolev embedding inequality.
	\begin{lem} \cite{cho2009sobolev}
		Let $d \geq 2$ and  $1 / 2 \leq \varsigma<1,$ then there exists $C$ such that for all radial $u \in H^{1}$
		\begin{equation} \label{radials}
			\sup _{x \in \mathbb{R}^{d} \backslash\{0\}}|x|^{d / 2-\varsigma}|u(x)| \leq C(d, \varsigma)\|u\|_{L^{2}( \mathbb{R}^d)}^{1-\varsigma}\|\nabla u\|_{L^{2}( \mathbb{R}^d)}^{\varsigma}.
		\end{equation}
	\end{lem}
	Later in this paper, we always  set $ \varsigma=1/2  $. By Lemma \ref{le4.1} and \eqref{radials},  we have
	\begin{align} \label{III}
		\begin{aligned}
			\textmd{III} &=-\frac{1}{2} \int_{\mathbb{R}^{d}}(-\Delta)^{\frac{1}{2}}\left(\Delta a_{R}(-\Delta)^{-\frac{1}{2}}|u|^{\alpha-1} u\right) u d x \\
			&=-\frac{1}{2} \int_{\mathbb{R}^{d}} \Delta a_{R}|u|^{\alpha+1} d x+\frac{1}{2} \int_{\mathbb{R}^{d}}\left[(-\Delta)^{\frac{1}{2}}, \Delta a_{R}\right]\left((-\Delta)^{-\frac{1}{2}}|u|^{\alpha-1} u\right) u d x \\
			& \geq-d \int_{|x| \leq \frac{R}{2}}|u|^{\alpha+1} d x-C \int_{|x| \geq \frac{R}{2}}|u|^{\alpha+1} d x\\
			&\qquad-C\left\|\nabla \Delta a_{R}\right\|_{L^{\infty}}\left\|(-\Delta)^{-\frac{1}{2}}|u|^{\alpha-1} u\right\|_{L^{2}}\|u\|_{L^{2}}\\
			&		\geq-d \int_{\mathbb{R}^{d}}|u|^{\alpha+1} d x-C\|u\|_{L^{\infty}\left(|x| \geq \frac{R}{2}\right)}^{\alpha-1}\|u\|_{L^{2}}^{2}-\frac{C}{R}\|u\|_{H^{1}}^{\alpha+1} \\
			&	\geq-d \int_{\mathbb{R}^{d}}|u|^{\alpha+1} d x-C\left(R^{-\frac{d-1}{2} (\alpha-1)}+R^{-1}\right).
		\end{aligned}
	\end{align}
	Similarly,
	\begin{equation} \label{II}
		\begin{aligned}
			\textmd{II} 	=&-2 \int_{\mathbb{R}^{d}} \partial_{j} \partial_{k}(-\Delta)^{-1}|u|^{\alpha-1} u \partial_{j} \partial_{k} a_{R} u d x \\
			&	-2 \int_{\mathbb{R}^{d}}\left[\partial_{j} \partial_{k}(-\Delta)^{-\frac{1}{2}}, \partial_{j} \partial_{k} a_{R}\right]\left((-\Delta)^{-\frac{1}{2}}|u|^{\alpha-1} u\right) u d x \\
			\geq & 2d \int_{|x| \leq \frac{R}{2}}|u|^{\alpha+1} d x-C \int_{|x| \geq \frac{R}{2}}|u |\left| \partial_{j} \partial_{k}(-\Delta)^{-1}\left| u\right|^{\alpha-1} u \right| d x \\	
			&	-C\left\|\nabla \partial_{j} \partial_{k} a_{R}\right\|_{L^{\infty}}\|u\|_{L^{2}}\left\|(-\Delta)^{-\frac{1}{2}}|u|^{\alpha-1} u\right\|_{L^{2}}\\
			\geq & 2 d \int_{\mathbb{R}^{d}}|u|^{\alpha+1} d x-C\|u\|_{L^{(2-l)^{\prime}}\left(|x| \geq \frac{R}{2}\right)}\left\|\partial_{j} \partial_{k}(-\Delta)^{-1}|u|^{\alpha-1} u\right\|_{L^{2-l}}-\frac{C}{R} \\
			\geq & 2d \int_{\mathbb{R}^{d}}|u|^{\alpha+1} d x-C\left(R^{-\frac{l(d-1)}{2(2-l)} }+R^{-1}\right),
		\end{aligned}
	\end{equation}
	where $ l=4/(d+4), (d-1)l/2(2-l)= \frac{(d+2)(d-1)}{d(d+4)} $. 
	
	Noticing that $ \left[(-\Delta)^{\frac{1}{2}}, x_{j}\right]=\mathscr{F}^{-1} \xi_{j} /(i|\xi|) \mathscr{F}=R_{j} $, one has
	\begin{equation}
		\begin{aligned}
			\textmd{I} =&-\int_{\mathbb{R}^{d}}(-\Delta)^{\frac{1}{2}}\left(2 x_{j} \partial_{j}(-\Delta)^{-\frac{1}{2}}|u|^{\alpha-1} u\right) u d x \\
			& \quad-\int_{\mathbb{R}^{d}}(-\Delta)^{\frac{1}{2}}\left(\left(\partial_{j} a_{R}-2 x_{j}\right) \partial_{j}(-\Delta)^{-\frac{1}{2}}|u|^{\alpha-1} u\right) u d x \\
			=&-2 \int_{\mathbb{R}^{d}} x_{j} u \partial_{j}|u|^{\alpha-1} u d x \\
			&-2 \int_{\mathbb{R}^{d}} u R_{j} \partial_{j}(-\Delta)^{-\frac{1}{2}}|u|^{\alpha-1} u d x \\
			&-\int_{\mathbb{R}^{d}}(-\Delta)^{\frac{1}{2}}\left(\left(\partial_{j} a_{R}-2 x_{j}\right) \partial_{j}(-\Delta)^{-\frac{1}{2}}|u|^{\alpha-1} u\right) u \chi_{R} d x\\
			&	-\int_{\mathbb{R}^{d}}(-\Delta)^{\frac{1}{2}}\left(\left(\partial_{j} a_{R}-2 x_{j}\right) \partial_{j}(-\Delta)^{-\frac{1}{2}}|u|^{\alpha-1} u\right) u (1-\chi_{R}) d x\\
			:&=\textmd{I}_1+\textmd{I}_2+\textmd{I}_3+\textmd{I}_4.
		\end{aligned}
	\end{equation}
	Directly calculation gives
	\begin{align} \label{II2}
		\textmd{I}_1=\frac{2\alpha d}{\alpha+1} \int_{\mathbb{R}^{d}}|u|^{\alpha+1} d x, \quad \textmd{I}_2=-2 \int_{\mathbb{R}^{d}}|u|^{\alpha+1} d x,
	\end{align}
	and
	\begin{equation} \label{I4}
		\begin{aligned}
			\textmd{I}_4=&-\int_{\mathbb{R}^{d}}\left(1-\chi_{R}\right) u\left(\partial_{j} a_{R}-2 x_{j}\right) \partial_{j}|u|^{\alpha-1} u d x \\
			&-\int_{\mathbb{R}^{d}}\left(1-\chi_{R}\right) u\left[(-\Delta)^{\frac{1}{2}},\left(\partial_{j} a_{R}-2 x_{j}\right)\right] \partial_{j}(-\Delta)^{-\frac{1}{2}}|u|^{\alpha-1} u d x \\
			=& \frac{\alpha}{\alpha+1} \int_{\mathbb{R}^{d}} \partial_{j}\left(\left(1-\chi_{R}\right)\left(\partial_{j} a_{R}-2 x_{j}\right)\right)|u|^{\alpha+1} d x\\
			&-\int_{\mathbb{R}^{d}}\left(1-\chi_{R}\right) u\left[(-\Delta)^{\frac{1}{2}},\left(\partial_{j} a_{R}-2 x_{j}\right)\right] \partial_{j}(-\Delta)^{-\frac{1}{2}}|u|^{\alpha-1} u d x \\
			\geq & C \int_{|x| \geq \frac{R}{2}}|u|^{\alpha+1} d x-C\left\|\left(1-\chi_{R}\right) u\right\|_{L^{(2-l)^\prime}}\left\|\partial_{j}(-\Delta)^{-\frac{1}{2}}\left(|u|^{\alpha-1} u\right)\right\|_{L^{2-l}}\\
			\geq & -C\left(R^{-\frac{(d-1) (\alpha-1)}{2}}+R^{-\frac{d-1}{2} \frac{l}{2-l}}\right),
		\end{aligned}
	\end{equation}
	and
	\begin{equation}
		\begin{aligned}
			\textmd{I}_3=&-\int_{\mathbb{R}^{d}}(-\Delta)^{\frac{1}{2}}\left(\chi_{R} u \right)\left(\partial_{j} a_{R}-2 x_{j}\right) \partial_{j}(-\Delta)^{-\frac{1}{2}}|u|^{\alpha-1} ud x\\
			\geq& -C\left\||x|(-\Delta)^{\frac{1}{2}}\left(\chi_{R} u\right)\right\|_{L^{(2-l)^{\prime}}(|x| \geq R / 2)}\left\|\partial_{j}(-\Delta)^{-\frac{1}{2}}\left(|u|^{\alpha-1} u\right)\right\|_{L^{2-l}}.
		\end{aligned}
	\end{equation}
	 Given $ x \in \mathbb{R}^{d},|x| \geq R / 2 $, we have
	\begin{equation}
		\begin{aligned}
			(-\Delta)^{\frac{1}{2}}\left(\chi_{R} u\right)(x)=& \frac{1}{(2 \pi)^{\frac{d}{2}}} \int_{\mathbb{R}^{d}}|\xi| \widehat{\chi_Ru}(\xi) e^{i x \cdot \xi} d \xi \\
			=& \frac{1}{(2 \pi)^{\frac{d}{2}}} \int_{\mathbb{R}^{d}} \chi(M \xi)|\xi| \widehat{\chi_{R} u}(\xi) e^{i x \cdot \xi} d \xi \\
			&+\frac{1}{(2 \pi)^{\frac{d}{2}}} \int_{\mathbb{R}^{d}}(1-\chi(M \xi))|\xi| \widehat{\chi_{R} u}(\xi) e^{i x \cdot \xi} d \xi \\
			:=& A_{1}+A_{2}.
		\end{aligned}
	\end{equation} 
	$  A_{1} $ can be estimated directly 
	\begin{align*}
		\left|A_{1}\right| \lesssim M^{-1} M^{-\frac{d}{2}}\left\|\chi_{R} u\right\|_{L^{2}} \lesssim M^{-1-\frac{d}{2}}.
	\end{align*}
	Integrating by parts, we obtain
	\begin{align*}
		\left|A_{2}\right| &=\left|\frac{1}{(2 \pi)^{d}} \int_{\mathbb{R}^{d}}(1-\chi(M \xi))|\xi| \int_{\mathbb{R}^{d}} \chi_{R} u(y) e^{i(x-y) \cdot \xi} d y d \xi \right|\\
		&=\left|(-1)^{N} \frac{1}{(2 \pi)^{d}} \int_{\mathbb{R}^{d}} \int_{\mathbb{R}^{d}} \Delta^{N}((1-\chi(M \xi))|\xi|) \frac{\chi_{R}(y) u(y)}{|x-y|^{2 N}} e^{i(x-y)}  d y d \xi\right|\\
		&\lesssim|x|^{-2 N} M^{2 N-d-1}\left\|\chi_{R} u\right\|_{L^{1}} \lesssim|x|^{-2 N} M^{2 N-d-1} R^{\frac{d}{2}},
	\end{align*}
	where $ N\in \mathbb{N}^\ast $.
	Let $ M=|x|^{2 N /(2 N-d / 2)} R^{-d /(4 N-d)} $, then
	\begin{equation}
		\begin{aligned}
			&\left\||x|(-\Delta)^{\frac{1}{2}}\left(\chi_{R} u\right)\right\|_{L(2-l)^{\prime}(|x| \geq R / 2)} \\
			&	\lesssim\left\||x||x|^{-2 N(d+2) /(4 N-d)} R^{d(d+2) /(8 N-2 d)}\right\|_{L^{(2-l)^{\prime}}(|x| \geq R / 2)} \\
			&\lesssim R^{-(\frac{1}{2-l}-\frac12)d},
		\end{aligned}
	\end{equation}
	where $ l=4/(d+4) $.
	Thus, 
	\begin{align} \label{I3}
		\textmd{I}_3 \gtrsim -R^{-\frac{d+2}{d+4}}.
	\end{align}

	Overall, we have obtained
	\begin{equation} \label{tsest1}
		\begin{aligned}
			\int_{\mathbb{R}^{d}} &|u|^{\alpha-1}  u(-\Delta)^{\frac{1}{2}}\left(\nabla a_{R} \cdot \nabla(-\Delta)^{-\frac{1}{2}} u+\frac{1}{2} \Delta a_{R}(-\Delta)^{-\frac{1}{2}} u\right) d x \\
			& \geq (-d+2d+\frac{2\alpha d}{\alpha+1}-2) \int_{\mathbb{R}^{d}}|u|^{\alpha+1} d x\\
			&-\left(\frac1R+R^{-\frac{d-1}{2} (\alpha-1)}+R^{-\frac{d-1}{2} \frac{l}{2-l}}+R^{-(\frac{1}{2-l}-\frac12)d}\right)\\
			& \geq  (d+\frac{2\alpha d}{\alpha+1}-2) \int_{\mathbb{R}^{d}}|u|^{\alpha+1} d x-C R^{-\theta}, 
		\end{aligned}
	\end{equation}
	where
	\begin{align*}
		\theta :=\min \left\{1,\frac{d-1}{2} (\alpha-1), \frac{(d+2)(d-1)}{d(d+4)}\right\}.
	\end{align*}
	This is \eqref{goal2}.
	\begin{flushright}
		$\square$
	\end{flushright}
	\begin{cor}  \label{cor}
		For any $ n\in \mathbb{Z}^+ $, we can find $ T_n \geq n $ such that 
		\begin{align} \label{tsest2}
			\|u\|_{L^{\alpha+1}([T_n,n+T_n],L^{\alpha+1})} \leq \frac{1}{2^n} .
		\end{align}
	\end{cor}
	\noindent{\bf{Proof:}}
	
	If not, then for some $ n_0\in \mathbb{Z}^+ $ and case $ d\geq 4 $, one has 
	\begin{align*}
		\|u\|_{L^{\alpha+1}([kn_0,(k+1)n_0],L^{\alpha+1})} > \frac{1}{2^{n_0}}, ~~\forall ~k\geq 1.
	\end{align*}
	For any $ K\in  \mathbb{Z}^+  $, it follows \eqref{tsinte1} that 
	\begin{align*}
		\frac{K}{2^{(\alpha+1)n_{0}}}<\int_{n_{0}}^{(K+1) n_{0}} \int_{\mathbb{R}^{d}}|u(t, x)|^{\alpha+1} d x d t \leq C K^{\frac{1}{1+\theta}} n_{0}^{\frac{1}{2+\theta}}
	\end{align*}
	which is a contradiction when $ K $ is sufficiently large. Similar discussion can show the results for case $ d=3 $.
	\begin{flushright}
		$\square$
	\end{flushright}	
	
	Now, we can be able to  obtain the space-time bound.
	\begin{prop}
		The $ H^1 $ wellposed solutions to \eqref{ch4_s1} satisfy 
		\begin{align}
			\left\|\left\langle\nabla\right\rangle^{s_c} v\right\|_{L_{t}^{q}\left(\mathbb{R}; L_{x}^{q}(\mathbb{R}^d)\right)}  \lesssim_{\left\|v_0\right\|_{H^1}} 1, \label{globalest}
		\end{align}
		where $  s_c=d / 2-2 / (\alpha-1), q=2(d+2) / d, $  $ 4/d+1\leq\alpha<(d+2)/(d-2) $ and $ d\geq 3 $.
	\end{prop}	
	
	\noindent{\bf{Proof:}}
	
	\eqref{ch4_inte1} can be written as 
	\begin{equation}
		\begin{aligned}			
			v(t)
			=& e^{-it \mathfrak{B}} v_{0} -i\int_0^{T_n}	e^{-i(t-\tau)\mathfrak{B}} \mathfrak{M}|\textmd{Re}(v)|^{\alpha-1}\textmd{Re}(v)d \tau\\
			&-i\int_{T_n}^{n+T_n} e^{-i(t-\tau)\mathfrak{B}} \mathfrak{M}|\textmd{Re}(v)|^{\alpha-1}\textmd{Re}(v)d \tau\\
			&-i\int_{n+T_n}^t	e^{-i(t-\tau)\mathfrak{B}} \mathfrak{M}|\textmd{Re}(v)|^{\alpha-1}\textmd{Re}(v)d\tau\\
			:=& 	e^{-it \mathfrak{B}}  v_{0}+D_1+D_2+D_3.
		\end{aligned}
	\end{equation}
	Denote $ \|v\|_{X(I)}=\left\|\left\langle\nabla\right\rangle^{s_c} v\right\|_{L_{t}^{q}\left(I; L_{x}^{q}\right)}, $  where $ 0\leq s_c <1 $. For any $ T>n+T_n $, it follows Strichartz estimates \eqref{stype4}, \eqref{stype5} that
	\begin{equation} \label{D3}
		\begin{aligned}
			\left\|D_{3}\right\|_{X\left([n+T_{n}, T]\right)} & \lesssim\left\|\left\langle\nabla\right\rangle^{s_c}|v|^{\alpha-1} v\right\|_{L^{p_{2}^{\prime}}\left([n+T_{n}, T]; L_{x}^{r_{2}^{\prime}}\right)} \\
			& \lesssim\left\|\left\langle\nabla\right\rangle^{s_c} v\right\|_{L^{q}\left(\left[n+T_{n}, T\right]; L_{x}^{q}\right)}\|v\|_{L^{q}\left(\left[n+T_{n}, T\right]; L_{x}^{r}\right)}^{\alpha-1} \\
			& \lesssim\|v\|_{X\left(\left[n+T_{n}, T\right]\right)}^{\alpha} \lesssim\|v\|_{X\left(\left[n+T_{n}, T\right]\right)}^{\alpha},
		\end{aligned}
	\end{equation}
	where $  p_{2}^{\prime}=q /\alpha, r_{2}^{\prime}=2 d(d+2) /\left(d^{2}+6 d+8-2 \alpha d\right), r=(\alpha -1)d(d+2) /(4+ 
	(3-\alpha) d).$
	
	Noticing that $ D_{2}=	e^{-it \mathfrak{B}} \left[	e^{i(n+T_n) \mathfrak{B}}  v\left(n+T_{n}\right)-	e^{iT_n \mathfrak{B}}  v\left(T_{n}\right)\right] $, one can use Strichartz estimates \eqref{stype4} and conservation law \eqref{ch4_est1} to get
	\begin{align}
		\left\|\left\langle\nabla\right\rangle D_{2}\right\|_{L^{q} ([n+T_n,T];L_{x}^{q})} \lesssim \left\|v\left(n+T_{n}\right)\right\|_{H^1}+\left\|v\left(T_{n}\right)\right\|_{H^1} \lesssim \mathcal{E}(0) \label{H1bounded}
	\end{align}
	and the Gagliardo–Nirenberg inequality, \eqref{tsest2}, Strichartz estimates \eqref{stype4}, \eqref{stype5}	to give
	\begin{equation} \label{D2}
		\begin{aligned}
			\left\|D_{2}\right\|_{X\left(\left[n+T_{n},T\right]\right)} & \lesssim\left\|\left\langle\nabla\right\rangle D_{2}\right\|_{L^{q} ([n+T_n,T];L_{x}^{q})}^{s_{\alpha}}\left\|D_{2}\right\|_{L^{q}\left(\left[n+T_{n},T\right]; L_{x}^{q}\right)}^{1-s_{\alpha}} \\
			& 	\lesssim\left\||u|^{\alpha-1} u\right\|_{L^{2}([T_n, n+T_{n}]; L_{x}^{\frac{2 d}{d+2}})}^{1-s_{\alpha}} \\
			&	\lesssim\|u\|_{L^{1+\alpha}\left(\left[T_n,n+T_{n}\right]; L_{x}^{1+\alpha}\right)}^{\frac{(1+\alpha) (1-s_{\alpha})}{2}}\|u\|_{L_t^{\infty} L_{x}^{\frac{(\alpha-1)d}{2}}}^{\frac{(\alpha-1)(1-s_{\alpha})}{2}} \\
			&	\lesssim\|u\|_{L^{1+\alpha}\left(\left[T_{n}, n+T_{n}\right]; L_{x}^{1+\alpha}\right)}^{\frac{(1+\alpha) (1-s_{\alpha})}{2}} \lesssim 2^{-n \frac{\left(1-s_{\alpha}\right)(1+\alpha)}{2}}.
		\end{aligned}
	\end{equation}
	
	To dispose $ D_1 $, we need the radial Strichartz estimates. For radial symmetry function $ \phi  \in L^2 (\mathbb{R}^d), d\geq 2$, the (3.16) in Corollary 3.4. of Z. Guo et al.  \cite{guo2018scattering} gives 
	\begin{align} \label{D15}
		\left\|e^{-i t \mathfrak{B} }  \phi\right\|_{L_{t}^{q} L_{x}^{r}} \lesssim\|\phi\|_{L_{x}^{2}}, ~d \geq 2
	\end{align}
	where $ (r,q) $  satisfy $ q\left(1/2-1/r\right)=1/(d-1),~ 2/(2d-1)\leq q\left(1/2-1/r\right)$,
	and it follows that
	\begin{align} \label{D14}
		\left\|D_{1}\right\|_{L^{2} (\left[n+T_{n},T\right]; L_{x}^{\frac{2 d-2}{d-2}})} \lesssim\left\|D_1(0)\right\|_{L^{2}} \lesssim 1.
	\end{align}
	Assume $ (2-\alpha) / 2 \leq 1 / r<(d-2) / 2 d, 2<r<\infty $, then 
	\begin{equation} \label{D13}
		\begin{aligned}
			\left\|D_{1}\right\|_{L_{t}^{\infty}\left(\left[n+T_{n},T\right]; L_{x}^{r}\right)} & \lesssim \int_{0}^{T_{n}}\left(T_{n}+n-s\right)^{-d\left(\frac{1}{2}-\frac{1}{r}\right)}\|u(s)\|_{L^{\alpha r^{\prime}}}^{\alpha} d s \\
			& \lesssim n^{1-d\left(\frac{1}{2}-\frac{1}{r}\right)},
		\end{aligned}
	\end{equation}	
	where we used $ \left\|e^{-i t\mathfrak{B}} f\right\|_{L^r} \leq c|t|^{-d\left(1 /2 -1 / r\right)}\|f\|_{L^{r^{\prime}}} $.
	
	If $r \geq 4-4 / d$, it follows interpolation inequality, H\"older inequality and Strichartz estimates \eqref{stype4}, \eqref{stype5} that
	\begin{equation} \label{D11}
		\begin{aligned}
			\left\|D_{1}\right\|_{L^{q}\left(\left[n+T_{n},T\right]; L_{x}^{q}\right)}
			& \lesssim\left\|D_{1}\right\|_{L_{t}^{2} (\left[n+T_{n},T\right];  L_{x}^{\frac{2d-2}{d-2}})}^{1-\theta_{1}-\theta_{2}}\left\|D_{1}\right\|_{L_{t}^{\infty}\left(\left[n+T_{n},T\right]; L_{x}^{r}\right)}^{\theta_{1}}\left\|D_{1}\right\|_{L_{t}^{\infty}( \left[n+T_{n},T\right];  L_{x}^{2})}^{\theta_{2}} \\
			& \lesssim n^{\theta_{1}-d \theta_{1}\left(\frac{1}{2}-\frac{1}{r}\right)},
		\end{aligned}
	\end{equation}
	where $\theta_{1}=(d-2) r /((d-1)(d+2)(r-2)), \theta_{2}=(d r-4 d+4) /((d-1)(d+$
	2) $(r-2))$, $ 1-\theta_1-\theta_{2}=d/(d+2) $. 
	If $r \leq 4-4 / d$, one has
	\begin{equation} \label{D12}
		\begin{aligned}
			\left\|D_{1}\right\|_{L_{t}^{q}\left(\left[0,T_{n}\right];L_{x}^{q}\right)} & \lesssim\left\|D_{1}\right\|_{L^{2} ([n+T_n,T]; L_{x}^{\frac{2d-2}{d-2}})}^{1-\theta_{1}-\theta_{2}}\left\|D_{1}\right\|_{L_{t}^{\infty}\left(\left[n+T_{n}, T\right]; L_{x}^{r}\right)}^{\theta_{1}}\left\|D_{1}\right\|_{L^{2}([n+T_n,T]; L_{x}^{\frac{2 d}{d-2}})}^{\theta_{2}} \\
			& \lesssim n^{\theta_{1}-d \theta_{1}\left(\frac{1}{2}-\frac{1}{r}\right)} ,
		\end{aligned}
	\end{equation}
	where $ \theta_{1}=2 /(d+2), \theta_{2} =d(4 d-r d-4) /((d-2)(d+2) r)  $. Similarly as \eqref{H1bounded}, one has $ 	\left\|\left\langle\nabla\right\rangle D_{1}\right\|_{L_{t}^{q} L_{x}^{q}} \lesssim \mathcal{E}(0) $, thus
	\begin{equation} \label{D1}
		\begin{aligned}
			\left\|D_{1}\right\|_{X\left(\left[n+T_{n},T\right]\right)} & \lesssim\left\|\left\langle\nabla\right\rangle D_{1}\right\|_{L_{t}^{q} L_{x}^{q}}^{s_{\alpha}}\left\|D_{1}\right\|_{L^{q}\left(\left[n+T_{n},T\right]; L_{x}^{q}\right)}^{1-s_{\alpha}} \\
			& \lesssim n^{\left(1-d\left(\frac{1}{2}-\frac{1}{r}\right)\right) \theta_{1}\left(1-s_{\alpha}\right)}
		\end{aligned}
	\end{equation}
	
	For any given $\varepsilon>0$, since $\left\|	e^{-it \mathfrak{B}}  v(0)\right\|_{X} \lesssim\|v(0)\|_{\dot{H}^{s_{\alpha}}}<\infty$, we can choose $n$
	sufficient large such that $\left\|	e^{-it \mathfrak{B}}  v_0\right\|_{X} \leq \varepsilon$ and follow \eqref{D2}, \eqref{D1} to get
	$$
	\left\|D_{1}\right\|_{X\left(\left[n+T_{n}, T\right]\right)}+\left\|D_{2}\right\|_{X\left(\left[n+T_{n}, T\right]\right)} \leq \varepsilon .
	$$
	Thus, we obtain
	$$
	\|v\|_{X\left(\left[n+T_{n}, T\right]\right)} \leq 2 \varepsilon+C\|v\|_{X\left(\left[n+T_{n}, T\right]\right)}^{1+\alpha}, \quad \forall T>0.
	$$
	By standard bootstrap argument, we arrive at  $v \in X(\mathbb{R})$.
	\begin{flushright}
		$\square$
	\end{flushright}
	
	\begin{cor}
		The $ H^1 $ wellposed solutions to \eqref{ch4_s1} satisfies 
		\begin{align} \label{globalest2}
			\left\| v\right\|_{S^1\left(\mathbb{R}\times\mathbb{R}^d\right)}  \lesssim_{ \left\|v_0\right\|_{H^1}} 1,
		\end{align}
		where $ S^1\left(\mathbb{R}\times\mathbb{R}^d\right) $ is defined as \eqref{S1n} and the $ \alpha $ in \eqref{ch4_s1} satisfies  $ 4/d+1\leq\alpha<(d+2)/(d-2) $ and $ d\geq 3 $.
	\end{cor}

	Combining with \eqref{globalest} and follow the standard analysis in Proposition 3.31 of T.Tao \cite{tao2006nonlinear}, we can obtain the \eqref{globalest2}, thus can prove Theorem \ref{th4}. In fact, for $ t_1<t_2 $, we have
	\begin{equation}
		\begin{aligned}
			\left\|	e^{it_1 \mathfrak{B}} v\left(t_{1}\right)-	e^{it_2 \mathfrak{B}} v\left(t_{2}\right)\right\|_{H^{1}} &\sim\left\|\int_{t_{1}}^{t_{2}} 	e^{i\tau \mathfrak{B}}  \left\langle\nabla\right\rangle \mathfrak{M}|\textmd{Re}(v)|^{\alpha-1}\textmd{Re}(v) d \tau\right\|_{L^{2}}\\
			&\lesssim\left\|\left\langle\nabla\right\rangle|v|^{\alpha-1} v\right\|_{L_{t}^{p_{2}^{\prime}}([t_1, t_2]; L_{x}^{r_{2}^{\prime}})} \\
			& \lesssim\left\|\left\langle\nabla\right\rangle v\right\|_{L_{t}^{q}\left(\left[t_1, t_2\right]; L_{x}^{q}\right)}\|v\|_{L_{t}^{q}\left(\left[t_1, t_2\right]; L_{x}^{r}\right)}^{\alpha-1} \\
			&\lesssim	\|v\|_{S^{1}\left([t_1,t_2] \times R^{d}\right)}^{\alpha}	\rightarrow 0, \text{     as   }  t_1,t_2 \rightarrow \infty,	
		\end{aligned}
	\end{equation}
	where $ q,p_2, r_2, r $ is defined as in \eqref{D3}. The rest is as those in Section \ref{s3.}.
	
	\begin{rem}
		One of the surprising points of Theorem \ref{th4} is that we can prove scattering for the mass critical  defocusing case $ \alpha=1+d/4, d\geq 3 $ with large radial initial data, while we do not  use of the  concentration compactness argument. This is different from those of defousing Schr\"odinger equation \eqref{nsl1}.	 By inequality \eqref{station}, we know that the low frequency of linear estimate of gBQ  is better than that of  Schr\"odinger equation \eqref{nsl1}. Then, we can get radial Strichartz estimate as \eqref{D15} ( see Guo et al. \cite{guo2018scattering}) and thus can interpolate between \eqref{D15} and energy estimate to obtain the smallness of $ D_1 $ in $ X\left(\left[n+T_{n}, T\right]\right) $. Note that Strichartz estimate as \eqref{D15} is impossible for Schr\"odinger equation \eqref{ch4_s1} within the analysis in Guo et al. \cite{guo2014improved,guo2018scattering}.
		
	\end{rem}
	
	\begin{appendices}
		\section{Analysis without transformation} \label{s5.}~
		
		We can prove Theorem \ref{th1}  without using the transform \eqref{trans0} in the Introduction. The Duhamel principle shows \eqref{ch4_1} is equivalent to the integral equation,
		\begin{equation} \label{ch4_inte3}
			u(t)=\cos\left(t\mathfrak{B}\right) u_{0}+\sin\left(t\mathfrak{B}\right)\mathfrak{B}^{-1} u_{1}-\beta\int_{0}^{t} \sin\left((t-s)\mathfrak{B}\right)\mathfrak{M}\left(|u|^{\alpha-1}u\right)d s,
		\end{equation}
		where  $  \mathfrak{B}, \mathfrak{M}, \cos\left(t\mathfrak{B}\right), \sin\left(t\mathfrak{B}\right) $ are defined as \eqref{ch4_trans1}, \eqref{cosB} respectively.
		As $ \cos\left(t\mathfrak{B}\right) f,$  $ \sin\left(t\mathfrak{B}\right) f $ are the real and imagine part of $  	e^{-it \mathfrak{B}} f$ respectively, the Strichartz estimates \eqref{stype4}, \eqref{stype5} may be still applied to $ \cos\left(t\mathfrak{B}\right), \sin\left(t\mathfrak{B}\right) $. One can obtain  local wellposedness results that are analogous to those in Section \ref{s2.} and \eqref{ch4_inte3} is equivalent to 
		\eqref{ch4_inte1}. For instance, by \eqref{stype4}, \eqref{stype5}, \eqref{nonlinear2}, we have 
		\begin{align*}
			\left\|u\right\|_{L^{p}([0,T];W^{1,q}(\mathbb{R}^{d}))}\lesssim
			\left\|u_0\right\|_{H^{1}(\mathbb{R}^{d})}&+
			\left\|\mathfrak{B}^{-1} u_1\right\|_{H^{1}(\mathbb{R}^{d})}+	\left\||u|^{\alpha}\right\|_{L^{p^\prime}([0,T];W^{1,q^\prime})}\nonumber\\
			\lesssim\left\|u_0\right\|_{H^{1}(\mathbb{R}^{d})}&+
			\left\|(-\Delta)^{-\frac12} u_1\right\|_{L^{2}(\mathbb{R}^{d})}+ \left\|\left\||u|^{\alpha-1}\right\|_{L^{\frac{d(\alpha+1)}{4}}}	\left\|u\right\|_{W^{1,q}}\right\|_{L^{p^\prime}[0,T]}\nonumber\\
			\lesssim\left\|u_0\right\|_{H^{1}(\mathbb{R}^{d})}&+
			\left\|(-\Delta)^{-\frac12} u_1\right\|_{L^{2}(\mathbb{R}^{d})}+	\left\|u\right\|_{L^{p}([0,T];W^{1,q})}^\alpha,
		\end{align*}

		A crucial tool to obtain these results is the stationary phase estimate derived in Gustafson et al. \cite{gustafson2006scattering}, Cho et al. \cite{cho2007small}
		\begin{align} \label{station}
			\sup _{x \in \mathbb{R}^{d}}\left|\int_{\mathbb{R}^{d}} e^{i\left(x \cdot \xi+t |\xi|\sqrt{1+|\xi|^2}\right)} \varphi\left(\frac{\xi}{N}\right) d \xi\right| \lesssim|t|^{-\frac{d}{2}}  \frac{N ^{\frac{d}{2}-1}}{\langle N\rangle^{\frac{d}{2}-1}},
		\end{align}
		where $ N \in 2^{\mathbb{Z}}, \varphi(\cdot ) $ is the Littlewood-Paley function. This immediately gives the the Strichartz estimates \eqref{stype1}, \eqref{stype2}, \eqref{stype3} and  improves the decay estimate introduced in
		\begin{align*}
			\sup _{x \in \mathbb{R}}\left|\int_{\mathbb{R}} e^{ i\left(x \cdot \xi+t |\xi|\sqrt{1+|\xi|^2}\right)} d \xi\right| \leqslant C\left(|t|^{-1 / 2}+|t|^{-1 / 3}\right),
		\end{align*}
		which is based on the  Van der Corput lemma, thus improves the Strichartz estimates in  Liu \cite{liu1997decay} and  the $ L_p$-$L_q $ estimates in Linares \cite{linares1993global}, which is a  special case of \eqref{stype4}, \eqref{stype5}  in one dimension. 
		Wang et al. \cite{wang2019cauchy} obtained the local  wellposedness for the generalized Boussinesq equation by adding damped terms $- \alpha\Delta n_t+\gamma \Delta^2 n_{tt},\alpha\geq0, \gamma>0, $  thus to give exponential decay term $ e^{-\frac{1}{2}|\xi|^{2}\left(\alpha+\gamma|\xi|^{2}\right) t}$ in the integral equation and obtain the  Strichartz type estimates. Here, we can get corresponding results of  \cite{wang2019cauchy} just within \eqref{station} and \eqref{stype1}-\eqref{stype3}.

		\section{Proof of Lemma \ref{le4.1} } \label{appb}~
		
		It is a corollary of Theorem 2 of Calder{\'o}n \cite{calderon1965commutators}, i.e.
		\begin{thm} \label{thB}
			Let $h(x)$ be homogeneous of degree $-d-1$ and locally integrable in $|x|>0 .$ Let $b(x)$ have first-order derivatives in $L^{r}, 1<r \leq \infty .$ Then, if $1<p<$ $\infty, 1<q<\infty, q^{-1}=p^{-1}+r^{-1}, h(x)$ is an even function and
			$$
			C_{\varepsilon}(f)=\int_{|x-y|>\varepsilon} h(x-y)[b(x)-b(y)] f(y) d y .
			$$
			$C_{\varepsilon}$ maps $L^{p}$ continuously into $L^{q}$ and $\|C(f)\|_{q} \leq c\|\operatorname{grad} b\|_{r}|| f \|_{p} \int \mid h(x) \mid d \nu$, where
			the integral is extended over $|x|=1, d \nu$ denotes the surface area of $|x|=1$, and c depends on $p$ and $r$ but not on $\varepsilon .$ Furthermore, as $\varepsilon$ tends to zero $C_{\varepsilon}(f)$ converges in norm in $L^{q}$.
		\end{thm} 
		
		In fact, 	denote $~c_{d, \alpha}:=\int_{\mathbb{R}^d}|x|^{-\alpha} e^{-|x|^{2} / 2} d x, $
		then
		\begin{align*}
			&	D=\sum_{j=1}^{d} R_{j} \partial_{j},~\mathscr{F}|\cdot|^{-1}=c_{1}|\cdot|^{-d+1} / c_{d-1},\\
			\mathscr{F}^{-1} \xi_{j} /(i|\xi|)&=-\partial_{j} \mathscr{F}^{-1}|\xi|^{-1}=-c_{1}(-d+1)|x|^{-d-1} x_{j} / c_{d-1}, \text{   for   } d\geq 3,\\
			R_{j} f(x)&=c \lim _{\varepsilon \rightarrow 0^{+}} \int_{B(0, \varepsilon)^{c}} \frac{y_{j} f(x-y)}{|y|^{d+1}} d y, \quad c=\frac{c_{1}(d-1)}{c_{d-1}} .
		\end{align*}
		For $f \in C_{0}^{\infty}( \mathbb{R}^d), \varphi \in C^{1, \alpha}( \mathbb{R}^d)$, one has
		\begin{align*}
			[D, \varphi] f&=D(\varphi f)-\varphi D f= \sum_{j=1}^{d} R_{j} \left(\partial_{j} \varphi  f\right)+R_{j}\left(\varphi \partial_{j} f\right)-\varphi R_{j} \partial_{j} f\\
			&=\sum_{j=1}^{d} R_{j} \left(\partial_{j} \varphi  f\right)+\sum_{j=1}^{d} \lim _{\varepsilon \rightarrow 0^{+}} \int_{B(0, \varepsilon)^{c}} \frac{y_{j}(\varphi(x-y)-\varphi(x)) \partial_{j} f(x-y)}{|y|^{d+1}} d y\\
			&=c \lim _{\varepsilon \rightarrow 0^{+}} \int_{B(0, \varepsilon)^{c}} \frac{-1}{|x|^{d+1}}(\varphi(x-y)-\varphi(x)) f(x-y) d y.
		\end{align*}
		Replace $ x-y $ with $ x $, we have 
		\begin{align}
			[D, \varphi] f(x)=c \lim _{\varepsilon \rightarrow 0^{+}} \int_{|x-y|>\varepsilon} \frac{1}{|x-y|^{d+1}}(\varphi(x)-\varphi(x-y)) f(y) d y.
		\end{align}
		Then, we can see \eqref{commu1} is the simplest case of Theorem \ref{thB}.
		
	\end{appendices}
	
	\begin{center}

	\end{center}

\end{document}